\documentclass[11pt,a4paper]{article}

\usepackage{amsmath}
\usepackage{amscd}
\usepackage{amssymb}
\usepackage{color}
\usepackage{latexsym}

\begin{document}

\newcommand{\CH}{{\mathcal H}}
\newcommand{\CT}{{\mathcal T}}
\newcommand{\CR}{{\mathcal R}}
\newcommand{\osigma}{\bar \sigma}
\newcommand{\ov}{\bar \bv}
\newcommand{\ovt}{\bar \vt}
\newcommand{\lpqt}{L_p(0,T;L_q(\Omega))}
\newcommand{\lpqtn}{L_p(0,T;L_q(\Omega)^{n-1})}
\newcommand{\BR}{\mathbb{R}}
\newcommand{\bG}{{\bold G}}
\newcommand{\bJ}{{\bold J}}

\newcommand{\qed}{\rightline{ $\square$}}
\renewcommand{\d}{\text{d}}
\renewcommand{\div}{\operatorname{div}}
\newcommand{\bv}{\boldsymbol{v}}
\newcommand{\bg}{\boldsymbol{g}}
\newcommand{\bn}{\boldsymbol{n}}
\newcommand{\tD}{\tilde{D}}
\newcommand{\btau}{\boldsymbol{\tau}}
\newcommand{\bT}{\boldsymbol{T}}
\newcommand{\bF}{\boldsymbol{F}}
\newcommand{\bH}{\boldsymbol{H}}
\newcommand{\bphi}{\boldsymbol{\varphi}}
\newcommand{\bI}{\boldsymbol{I}}
\newcommand{\bq}{\boldsymbol{q}}
\newcommand{\bw}{\boldsymbol{w}}
\newcommand{\bD}{\boldsymbol{D}}
\newcommand{\bA}{\boldsymbol{A}}
\newcommand{\bff}{\boldsymbol{f}}
\newcommand{\bomega}{\boldsymbol{\omega}}
\newcommand{\balpha}{\boldsymbol{\alpha}}
\newcommand{\bS}{\boldsymbol{S}}
\newcommand{\bb}{\boldsymbol{b}}
\newcommand{\bu}{\boldsymbol{u}}
\newcommand{\bnul}{\boldsymbol{0}}
\newcommand{\bpsi}{\boldsymbol{\psi}}
\newcommand{\bPhi}{\boldsymbol{\Phi}}
\newcommand{\de}{\partial}
\newcommand{\erz}{\mbox{e}^{z}}
\newcommand{\erzj}{\mbox{e}^{z^j}}
\newcommand{\erzjN}{\mbox{e}^{z^j_N}}
\newcommand{\etrzj}{\mbox{e}^{\widetilde{z^j}}}
\newcommand{\erzBj}{\mbox{e}^{Bz^j}}
\newcommand{\etrzBj}{\mbox{e}^{B\widetilde{z^j}}}
\newcommand{\emrzj}{\mbox{e}^{-z^j}}
\newcommand{\emrzjm}{\mbox{e}^{-z^{j-1}}}
\newcommand{\etmrzj}{\mbox{e}^{-\widetilde{z^j}}}
\newcommand{\etmrzjm}{\mbox{e}^{-\widetilde{z^{j-1}}}}
\newcommand{\erzjm}{\mbox{e}^{z^{j-1}}}
\newcommand{\erl}{\mbox{e}^{r_l}}
\newcommand{\erlN}{\mbox{e}^{r_{l,N}}}
\newcommand{\etrlN}{\mbox{e}^{\widetilde{r_{l,N}}}}
\newcommand{\erlj}{\mbox{e}^{r_l^j}}
\newcommand{\erk}{\mbox{e}^{r_k}}
\newcommand{\ezk}{\mbox{e}^{z_k}}
\newcommand{\erkj}{\mbox{e}^{r_k^j}}
\newcommand{\etrkj}{\mbox{e}^{\widetilde{r_k^j}}}
\newcommand{\etrkN}{\mbox{e}^{\widetilde{r_{k,N}}}}
\newcommand{\etrlj}{\mbox{e}^{\widetilde{r_l^j}}}
\newcommand{\erkjm}{\mbox{e}^{r_k^{j-1}}}
\newcommand{\erkN}{\mbox{e}^{r_{k,N}}}
\newcommand{\erkjN}{\mbox{e}^{r_{k,N}^j}}
\newcommand{\erkNd}{\mbox{e}^{r_{k,N,\delta}}}
\newcommand{\erki}[1]{\mbox{e}^{r_{k,#1}}}
\newcommand{\ezki}[1]{\mbox{e}^{z_{k,#1}}}

\newcommand{\Curl}{\operatorname{curl}}
\newcommand{\lap}{\Delta}
\newcommand{\rot}{\operatorname{rot}}
\newcommand{\Div}{\operatorname{div}}



\newcommand{\pf}{{\noindent\it Proof.~}}
\newcommand{\Ov}[1]{\overline{#1}}
\newcommand{\DC}{C^\infty_c}
\newcommand{\Un}[1]{\underline{#1}}
\newcommand{\bo}{| _{\partial\Omega}}
\newcommand{\vC}{{C}}
\newcommand{\vy}{Y_{k}}
\newcommand{\vF}{\boldsymbol{F}}
\newcommand{\bbf}{{\vc f}}
\newcommand{\vQ}{\vc{Q}}
\newcommand{\vS}{\boldsymbol{S}}
\newcommand{\vw}{\vc{w}}
\newcommand{\vv}{\vc{v}}
\newcommand{\hvw}{\hat{\omega}}
\newcommand{\vye}{\vy_\ep}
\newcommand{\Sd}{\{O,F,P\}}
\newcommand{\SD}{\{O,F,P,D\}}
\newcommand{\vr}{\varrho}
\newcommand{\tvr}{\tilde{\varrho}}
\newcommand{\hvr}{\hat{\varrho}}
\newcommand{\vp}{\varphi}
\newcommand{\tp}{\tilde{p}}
\newcommand{\ex}[1]{ \left< #1 \right>}
\newcommand{\vrn}{\vr_n}
\newcommand{\vte}{\vt_\ep}
\newcommand{\vun}{\vu_n}
\newcommand{\vrd}{\vr_\delta}
\newcommand{\vra}{\vr_{A}}
\newcommand{\vram}{\vr_{A-}}
\newcommand{\vrb}{\vr_{B}}
\newcommand{\tvra}{\tilde{\vr}_{A}}
\newcommand{\tvrb}{\tilde{\vr}_{B}}
\newcommand{\tvrn}{\tilde{\vr}_{n}}
\newcommand{\tg}{\tilde{g}}
\newcommand{\vtd}{\vt_\delta}
\newcommand{\vud}{\vu_\delta}
\newcommand{\vt}{\vartheta}
\newcommand{\tvt}{\tilde{\vartheta}}
\newcommand{\hvt}{\hat{\vartheta}}
\newcommand{\vu}{\boldsymbol{u}}
\newcommand{\ve}{{\cal E}}
\newcommand{\vd}{\boldsymbol{d}}
\newcommand{\vdk}{\boldsymbol{d}_k}
\newcommand{\vc}[1]{{\bf #1}}
\newcommand{\vcg}[1]{{\pmb #1}}
\newcommand{\Grad}{\nabla}
\newcommand{\pvr}{{{\partial} \over {\partial \vr}}}
\newcommand{\pvt}{{{\partial} \over {\partial \vt}}}
\newcommand{\pt}{\partial_{t}}
\newcommand{\ptb}[1]{\partial_{t}(#1)}
\newcommand{\Dt}{\frac{ d}{dt}}
\newcommand{\tn}[1]{\mbox {\F #1}}
\newcommand{\dx}{{\rm d} {x}}
\newcommand{\dt}{{\rm d} t }
\newcommand{\iit}{\it}
\newcommand{\Rm}{\mbox{\FF R}}
\newcommand{\dxdt}{\dx \ \dt}
\newcommand{\lr}[1]{\left( #1 \right)}
\newcommand{\intO}[1]{\int_{\Omega} #1 \ \dx}
\newcommand{\intOj}[1]{\int_{\Omega^{1}} #1 \ \dx}
\newcommand{\intOd}[1]{\int_{\Omega^{2}} #1 \ \dx}
\newcommand{\intOe}[1]{\int_{\Omega_\ep} #1 \ \dx}
\newcommand{\intOB}[1]{\int_{\Omega} \left( #1 \right) \ \dx}
\newcommand{\intRN}[1]{\int_{R^3} #1 \ \dx}
\newcommand{\intR}[1]{\int_{R} #1 \ \dt}
\newcommand{\intRR}[1]{\int_R \int_{R^3} #1 \ \dxdt}
\newcommand{\intT}[1]{\int_0^T #1 \ \dt}
\newcommand{\intTO}[1]{\int_0^T\!\!\!\! \int_{\Omega} #1 \ \dxdt}
\newcommand{\intt}[1]{\int_0^t #1 \ \dt}
\newcommand{\inttO}[1]{\int_0^t\!\!\!\! \int_{\Omega} #1 \ \dxdt}
\newcommand{\inttauO}[1]{\int_0^\tau\!\!\!\! \int_{\Omega} #1 \ \dxdt}
\newcommand{\intTOj}[1]{\int_0^T\!\!\!\! \int_{\Omega^{1}} #1 \ \dxdt}
\newcommand{\intTOd}[1]{\int_0^T\!\!\!\! \int_{\Omega^{2}} #1 \ \dxdt}
\newcommand{\intTOB}[1]{ \int_0^T\!\!\!\! \int_{\Omega} \left( #1 \right) \ \dxdt}
\newcommand{\inttOB}[1]{ \int_0^t\!\!\!\! \int_{\Omega} \left( #1 \right) \ \dxdt}
\newcommand{\sumkN}[1]{\sum_{k=1}^3 #1}
\newcommand{\sumlN}[1]{\sum_{l=1}^3 #1}
\newcommand{\blue}[1]{\textcolor{blue}{ #1}}
\newcommand{\red}[1]{\textcolor{red}{ #1}}
\newcommand{\eq}[1]{\begin{equation}
\begin{split}
#1
\end{split}
\end{equation}}
\newcommand{\eqh}[1]{\begin{equation*}
\begin{split}
#1
\end{split}
\end{equation*}}
\newcommand{\cB}{{\cal{B}}}
\newcommand{\cR}{{\cal{R}}}
\newcommand{\cF}{{\cal{F}}}
\newcommand{\cG}{{\cal{G}}}
\newcommand{\D}{{\vc{D}}}
\newcommand{\ep}{\varepsilon}
\newcommand{\dd}{\delta}
\newcommand{\R}{\mathbb{R}}
\newcommand{\C}{\mathbb{C}}
\newcommand{\vn}{{\bf n}}
\newcommand{\pd}{\partial}
\newcommand{\bV}{{\bf V}}
\newcommand{\bk}{{\bf k}}
\newcommand{\kv}{{\bk_{\bv}}}
\newcommand{\dv}{{\rm div}}
\newtheorem{thm}{Theorem}[section]
\newtheorem{lem}[thm]{Lemma}
\newtheorem{prop}[thm]{Proposition}
\newtheorem{df}{Definition}
\newtheorem{rmk}[thm]{Remark}

\title{On the isothermal compressible multi-component mixture flow:\\
 the local existence and maximal $L_p-L_q$ regularity of solutions}
\author{T. Piasecki\footnote{Institute of Applied Mathematics and Mechanics, University of Warsaw, 
ul. Banacha 2, 02-097 Warszawa, Poland. E-mail: {t.piasecki@mimuw.edu.pl}.
Supported by the Top Global University Project and the Polish National Science Centre grant 2018/29/B/ST1/00339.}, 
Y. Shibata\footnote{Department of Mathematics, Waseda University, Ohkubo 3-4-1, Shinjuku-ku, Tokyo 169-8555, Japan. 
Adjunct faculty member in the Department of Mechanical 
Engineering and Materias Science, University of Pittsburgh. E-mail: {yshibata@waseda.jp}.
Partially supported by  
JSPS Grant-in-aid for Scientific Research (A) 17H0109
and Top Global University Project.}, 
E. Zatorska\footnote{Department of Mathematics, University College London, Gower Street,  London WC1E 6BT, United Kingdom.  E-mail: {e.zatorska@ucl.ac.uk}.
Supported by the Top Global University Project and the Polish Government MNiSW research grant 2016-2019 "Iuventus Plus"  No.  0888/IP3/2016/74.}}

\date{}

\maketitle

\noindent{\bf{Abstract:}} We consider the initial-boundary value problem for the system of equations describing  the flow of compressible isothermal mixture of arbitrary large number of components. The system consists of the compressible Navier-Stokes equations and a subsystem of  diffusion equations for the species. The subsystems are coupled by the form of the pressure and the strong cross-diffusion effects in the diffusion fluxes of the species. Assuming the existence  of solutions to the symmetrized and linearized equations, proven in \cite{PSZ2}, we derive the estimates for the nonlinear equations and prove the local-in-time existence and maximal $L_p-L_q$ regularity of solutions.
\normalsize

\vspace{2cm}

\section{Introduction}
\subsection{Setting of the problem}
We consider the system of equations describing the motion of an isothermal mixture of compressible gases 
        \begin{equation} \label{1.1}
     \left.  
      \begin{array}{r}
            \vspace{0.2cm}
        \pt\vr+\Div (\vr \vu) = 0\\
            \vspace{0.2cm}
        \ptb{\vr\vu}+\Div (\vr \vu \otimes \vu) - \Div \bS+ \Grad p =\vc{0}\\
            \vspace{0.2cm}
        \pt{\vr_k}+\Div (\vr_{k} \vu)+ \Div \bF_k  =  0
            \end{array}\right\}\quad\mbox{in}\ (0,T)\times\Omega
    \end{equation}
in the regular domain $\Omega\subset \R^3$, supplied with boundary conditions 
\begin{equation} \label{bc}
\vu=0, \; \bF_k\cdot\bn=0 \quad\mbox{on}\ (0,T)\times\partial\Omega 
\end{equation}	
and initial condition 
\begin{equation} \label{bc}
\vu|_{t=0}=\vu^0, \quad \vr_k|_{t=0}=\vr_k^0,\; k=1 \ldots n \quad \mbox{in} \; \Omega.
\end{equation}
 Above, in system \eqref{1.1}, $\vr$ denotes the mass density of the mixture
\begin{equation}  \label{rho}
\vr=\sum_{k=1}^3\vr_k,
\end{equation}  
$\vu$ is the mean velocity of the mixture, and $\vr_k$ is the density of the $k$-th constituent.  
The remaining quantities:  
the stress tensor $\bS$, the total internal pressure $p$, and the diffusion fluxes $\bF_k$ 
are determined as functions of $(\vu,\vr,\vr_k)$ by constitutive relations which will be specified later.

The first equation of system \eqref{1.1}, usually called the continuity equation,  describes the balance of the mass, 
and the second equation expresses the balance of the momentum. 
The last $n$ equations describe the balances of masses  of separate constituents (species). 
Note that the system of equations cannot be independent, as the last $n$ equations must sum up to the continuity 
equation. Thus, here we meet a serious  mathematical obstacle,  the subsystem $(\ref{1.1})_4$ is 
degenerate parabolic in terms of $\vr_k$.

\medskip
{\bf The stress tensor.} 
The viscous part of the stress tensor obeys the {\it Newton rheological law}
	\begin{equation}\label{chF:Stokes}
	\vS(\vu)= 2\mu\D(\vu)+\nu\Div \vu{\bf{I}},
	\end{equation}
where $\D(\vu)=\frac{1}{2}\left(\Grad \vu+(\Grad \vu)^{T}\right)$ and the nonnegative viscosity coefficients. 

\medskip
{\bf Internal pressure.} 
The internal pressure of the mixture is determined through the Boyle law, when the temperature is constant it is given by
\begin{equation} \label{intpre} 
p(\vr_{1},\ldots,\vr_{n})=\sum_{k=1}^{n}p_{k}(\vr_{k})=\sum_{k=1}^{n}\frac{\vr_{k}}{m_{k}};
\end{equation}
above, $m_{k}$ is the molar mass of the species $k$, and for simplicity, we set the gaseous constant  equal to 1.


\medskip
{\bf Diffusion fluxes.} 
A key element of the presented model is the structure of laws governing cross-diffusion processes in the mixture. 
The diffusion fluxes are given explicitly in the form
	\begin{equation}\label{eq:diff}
	\bF_{k}=-\sum_{l=1}^{n} \vC_{kl}\vd_l, \quad k=1,...n,
	\end{equation}
where $\vC_{kl}$ are multicomponent flux diffusion coefficients and $\vdk=(d_{k}^{1},d_{k}^{2},d_{k}^{3})$ is the species $k$ diffusion force 
	\begin{equation}\label{eq:}
	d_{k}^{i}=\Grad_{x_{i}}\left({p_{k}\over p}\right)+\left({p_{k}\over p}-{\vr_{k}\over \vr}\right)\Grad_{x_{i}} \log{p}= \frac{1}{p}\lr{\Grad_{x_i} p_k-\frac{\vr_k}{\vr}\Grad_{x_i}p}.
	\end{equation}
Moreover, we assume that $\sumkN\bF_k=\vc{0}$, pointwisely.
The main properties of the flux diffusion matrix $C$ are
\begin{equation} \label{prop_C}
C{\cal Y}={\cal Y}C^{T},\quad
       N(C)=\mbox{lin}\{\vec Y\},\quad
      R(C)={U}^{\bot},
\end{equation}
where $Y_k=\frac{\vr_k}{\vr}$, ${\cal Y}=\mbox{diag}(Y_{1},\ldots,Y_{N})$, $\vec Y=(Y_1,\ldots,Y_n)^t$,
$\mbox{lin}\{\vec Y\}=\{t\vec Y\: \; t \in \R \}$,
$N(C)$ is the nullspace of $C$,  $R(C)$ is the range of $C$, 
$\vec U=(1,\ldots,1)^{T},$ and  ${U}^{\bot}$ is the orthogonal complement of $\mbox{lin}\{\vec U\}$.
The second property in \eqref{prop_C} implies 
$$
\sum_{l=1}^3 \frac{1}{p}C_{kl} \frac{\vr_k}{\vr}\nabla p=\frac{\nabla p}{p}\sum_{l=1}^3 C_{kl}Y_l=0,\quad k=1,\ldots, n,
$$
therefore \eqref{eq:diff}, \eqref{eq:} are reduced to
\begin{equation} \label{eq:diff1}
\bF_{k}=-\frac{1}{p}\sum_{l=1}^{n} \vC_{kl}\nabla p_l.
\end{equation}
We also define 
\begin{equation} \label{def_D}  D_{kl}=\frac{C_{kl}}{\vr Y_k},
\end{equation}
thus the properties of $C$ \eqref{prop_C}
imply 
\begin{equation} \label{prop_D}
 D= D^{T},\quad D \geq 0, \quad
       N(D)=\mbox{lin}\{\vec Y\},\quad
      R(D)={Y}^{\bot}.
\end{equation}
The first property results from $C_{kl}Y_l=C_{lk}Y_k$, the third from the fact that ${\cal Y}$
is diagonal. Next, $p \in R(\tilde D) \iff p_k=\frac{1}{Y_k}\sum_l C_{kl}q_l$ 
for some $q \in \mathbb{R}^3$. 
Finally $D$ is positive definite over $U^\bot$. 

\medskip
{\bf Exemplary diffusion matrix.} 
An example of matrix $C$ satisfying conditions \eqref{prop_C} that will be distinguished throughout the paper is
	\begin{equation}\label{Cform}
	\vC =\left(
		\begin{array}{cccc}
			Z_{1} & -Y_{1} & \ldots & -Y_{1}\\
			-Y_{2} & Z_{2}  & \ldots & -Y_{2}\\
			\vdots & \vdots & \ddots & \vdots \\
			-Y_{n} & -Y_{n}& \ldots & Z_{n}
		\end{array} \right),
	\end{equation}
where $Z_{k}=\sum_{{i=1} \atop {i\neq k}}^{n} Y_{i}$.\\
Using expressions for the diffusion forces \eqref{eq:diff1} and the properties of this matrix one can rewrite \eqref{eq:diff} into the following form
		\eq{\label{difp}
		\bF_k=-\frac{1}{p}\lr{\Grad p_k-Y_k\Grad p}.
		}
Clearly for $C$ given by \eqref{Cform}, the matrix $D_{kl}=\frac{C_{kl}}{\vr Y_k}$ is symmetric and positive semi-definite.

%

\subsection{Discussion of the known results}
The main result of this paper concerns the local well-posedness  of system \eqref{1.1} in the maximal $L_p-L_q$ regularity setting.
The local well-posedness as well as global well-posedness for small data  for two-species variant of system \eqref{1.1} have been  shown in authors' previous work  \cite{PSZ}. There 
the  so-called normal form, considered earlier e.g. in \cite{VG}, allows to immediately write a parabolic equation for one of the species densities. The aim of this paper is to generalize this result to the system with arbitrary number of constituents, however still isothermal. The key difference is that in the two species case the part corresponding to diffusion flux is reduced to a single parabolic equation, while now we obtain only a symmetrized system. Nevertheless, the properties of $D$ imply only nonnegativity of its leading order part so an important step is to show its parabolicity. Dealing with the systems of species instead of single equation also  requires serious modifications in the linear theory. 

The mathematical investigation of multicomponent flows dates back to analysis of a two component 
incompressible model assuming Fick law, hence no cross-diffusion, see among others \cite{BdV1}
for inviscid fluid and \cite{BdV2}-\cite{BdV3} in the viscous case.

In the previous results devoted to the complete mixture model, see Giovangigli and Massot \cite{GM1,GM2}, the local smooth solutions and global smooth solutions around constant equilibrium states were considered. Their method of proof was based on normal form of equations, hyperbolic-parabolic estimates  and on local strict dissipativity of linearized systems.
It can be seen as an application of more abstract theory proposed for the hyperbolic-parabolic systems of conservation laws by Kawashima and Shizuta \cite{K84,KS88}. 

When the species equations are decoupled from the fluid equations, the resulting system of PDEs is  related to the Stefan-Maxwell system analyzed for example in \cite{B2010, HMPW13}. In both of these papers the isobaric isothermal systems are considered with the barycentric velocity being equal to $0$. This means that, in comparison with the system of last $n$ equations from \eqref{1.1}, the convective term $\Div(\vr_k\bu)$ is absent and the variation of total pressure in the  diffusion fluxes \eqref{difp} is neglected. Essential difference between these systems is that in the present case the diffusion fluxes are explicit combination of diffusion deriving forces, while for the Stefan-Maxwell system the flux-forces relations need to be first inverted. This can be done using the Perron-Frobenius theory as first noticed in \cite{VG0}. With this at hand, the local-in-time well-posedness and maximal $L_p$ regularity follow from classical results of Amann \cite{Amann} or Pr\"uss \cite{Pruss}. In the approach presented in the present paper we rather relate on the alternative approach of the second author and collaborators \cite{ES1, SS2, Murata, MS16, S17, SS1} tailored to the compressible fluid systems. The main result of this paper is maximal $L_p-L_q$ regularity of solutions to \eqref{1.1}, but it relies on the proof of existence of relevant solutions to the linearized system.  The latter result is proved in our other article 
\cite{PSZ2} mostly for the sake of brevity, but also as it can be of independent interest. Indeed, it applies to whole class of symmetric parabolic systems satisfying certain regularity assumptions on the coefficients, therefore it is likely to be used in other contexts.

 As far as maximal $L_p-L_q$ regularity is concerned, the coupling between Stefan-Maxwell and the fluid equations, was so far considered only for the incompressible Navier-Stokes system, see \cite{BP2017}. It was also proven, independently in \cite{CJ13} and \cite{MT13}, that the incompressible Navier-Stokes-Stefan-Maxwell system possesses a global-in-time weak solution with arbitrary data .
The approach employed by Chen and J\"ungel in \cite{CJ13} relies on a certain symmetrization of the species subsystem with one of equations eliminated, see also \cite{JS13}. They have noticed that such reformulation allows to deduce parabolicity
in terms of the so-called entropic variables. See also \cite{Jungel} for an overview of different problems where a similar approach can be applied. 
The idea of our approach is similar, however the change of variables we propose is slightly different, in the spirit of normal variables from \cite{VG}.
Concerning analogous results for the compressible Navier-Stokes-Stefan-Maxwell system, the existence of weak solutions is so far known either for stationary flow of species with the same molar masses \cite{EZ, GPZ, PP1, PP2}, or for exemplary diffusion matrix $C$ and stress tensor $\vS$ with density-dependent viscosity coefficient \cite{EZ2,EZ3,MPZ1, MPZ2}. There are also relevant results for multi-component systems with diffusion fluxes in the form of the Fick law \cite{FPT}.

\subsection{Notation and functional spaces}
Let us summarize notation used in the paper. We use standard notation $H^k_p, \; k \in \mathbb{N}$ for Sobolev spaces. For a Banach space $X$, by $L_p(0,T;X)$ we denote a Bochner space and 
$$
H^1_p(0,T;X)=\{ f \in L_p(0,T;X): \; \de_t f \in L_p(0,T;X)\}.
$$
Furthermore, for $s \in \R$ a Bessel space $H^{s}_p(\R,X)$ is a space of $X$-valued functions for which 
\begin{equation*}
\|f\|_{H^{s}_p(\BR, X)}
= \Bigl(\int_\R \|\mathcal{F}^{-1}[(1+\tau^2)^{s/2}\mathcal{F}[f](\tau)]
\|^p\,{\rm d}\tau\Bigr)^{1/p} < \infty.      
\end{equation*}
We also recall that for $0<s<\infty$ and $m$ a smallest integer larger than $s$ we define Besov spaces on domains as intermediate spaces
\begin{equation} \label{def:bsqp0} 
B^{s}_{q,p}(\Omega)=(L_q(\Omega),H^m_q(\Omega))_{s/m,p},
\end{equation}
where $(\cdot,\cdot)_{s/m,p}$ is the real interpolation functor, see \cite[Chapter 7]{Ad}. In particular,
\begin{equation} \label{def:bsqp}
B^{2(1-1/p)}_{q,p}(\Omega)=(L_q(\Omega),H^2_q(\Omega))_{1-1/p,p}=(H^2_q(\Omega),L_q(\Omega))_{1/p,p}.    \end{equation}
Next, for abbreviation and clarity we introduce the following notation:
\begin{enumerate}
\item We will denote by $E(T)$ a continuous function of $T$ s.t. $E(0)=0$. Moreover, we use $C$ to denote a generic positive constant, or we use $C(X,Y)$ to specify the dependence of parameters $X$ and $Y$.
\item By $\vec\cdot$ we denote an $(n-1)$-vector  of functions, for example $\vec\vt=(\vt_1,\ldots,\vt_{n-1})^\top$.
\item We introduce the norms describing regularity of our solutions; for $T>0$ we define:
\eq{ \label{def:norm}
[{\bv}]_{T,1}&:=\|{\bv}\|_{L_p(0,T;H^2_q(\Omega))}+\|\de_t {\bv}\|_{\lpqt}, \\
[ \sigma]_{T,2}&:=\|\sigma\|_{H^1_p(0,T;H^1_q(\Omega))},\\
[\sigma,{\bv},{\vec\vt}]_T&:=[{\bv}]_{T,1}+[\sigma]_{T,2}+\sum_{k=1}^{n-1}[\vt_k]_{T,1}.
}
Then, for given $T,M>0$ we define the sets in the functional spaces:
\begin{align} \label{def:H12}
{\mathcal H}_{T,M}^1=\{ {\bv}: [{\bv}]_{T,1}\leq M \}, \qquad
{\mathcal H}_{T,M}^{2}=\{\sigma: \; [\sigma]_{T,2} \leq M\}
\end{align}
and
\begin{equation}\label{def:H}
{\mathcal H}_{T,M} = \left\{(\sigma, {\bv}, \vec\vt):  
\quad (\sigma, {\bv}, \vt_k)|_{t=0} = (0, \bu^0, h^0_k) \quad\text{in $\Omega$}, \quad 
[\sigma,{\bv},\vec\vt]_T \leq M\right\}.
\end{equation}
\end{enumerate}

\section{Symmetrization and main result}
The main result of this paper is the  the local well-posedness   in the maximal $L_p-L_q$ regularity setting of certain reformulation of system \eqref{1.1} \eqref{sys:normal}.
This reformulation  is similar to the normal form derived in (\cite{VG}, Chapter 8) for the complete system 
with thermal effects. In case of constant temperature derivation of the symmetrized equations  can be simplified considerably, and,  to make our paper self contained, we show in the Appendix the following result
\begin{prop} \label{thm:main0}
Let $(\vr,\vu,\vr_1,\ldots,\vr_n)$ be a regular solution to system (\ref{1.1}-\ref{rho}) such that
\eq{\label{all_pos}
\{\vr_1>C,\ldots,\vr_n>C\}}
for some constant $C>0$. Then the change of unknowns
\begin{equation} \label{def:psi}
(\vr,h_1,\ldots,h_{n-1})=\lr{
\sum_{i=1}^3 \vr_1,\log\lr{\frac{\vr_2^{\frac{1}{m_2}}}{\vr_1^{\frac{1}{m_1}}}},\ldots,\log\lr{\frac{\vr_n^{\frac{1}{m_n}}}{\vr_1^{\frac{1}{m_1}}}
}}=:\Psi(\vr_1,\ldots,\vr_n).
\end{equation}
is a diffeomorphism, and the system \eqref{1.1} is transformed to 
\eq{ \label{sys:normal}
&\pt\vr+\Div(\vr\vu)=0,\\
&\vr \de_t \vu+\frac{\vr\nabla\vr}{\Sigma_\vr}+\sum_{l=2}^3\lr{\vr_l-\frac{m_l\vr_l\vr}{\Sigma_\vr}}\nabla h_{l-1}+\vr(\vu\cdot \nabla)\vu
=\mu \Delta \vu + (\mu+\nu)\nabla\dv\vu,\\
&\sum_{l=1}^{n-1} \cR_{kl}(\de_t h_l+ \vu \cdot \nabla h_l) + \lr{\vr_{k+1}-\frac{m_{k+1}\vr_{k+1}\vr}{\Sigma_\vr}}\Div \vu
=\Div \left( \sum_{l=1}^{n-1}\cB_{kl}\nabla h_l\right),
}
with the boundary conditions
\begin{equation} \label{bc:normal}
\bu=0, \quad \sum_{l=1}^{n-1}\cB_{kl}\nabla h_l \cdot \bn = 0, \quad k=1,\ldots,n-1,\quad\mbox{on}\ (0,T)\times\partial\Omega,
\end{equation}
and the initial conditions
\begin{equation} \label{ic:normal}
(\bu,\vr,\{h_k\}_{k=1,\ldots,n-1})|_{t=0}=(\bu^0,\vr^0,\{h_k^0\}_{k=1,\ldots,n-1}) = \Psi(\vr_{1}^0(x),\ldots \vr_{n}^0(x)),
\end{equation}
where
\begin{equation} \label{def:sigma}
\Sigma_\vr=\sum_{k=1}^3 m_k \vr_k
\end{equation}
and  $\cR$ and $\cB$ are $(n-1)\times(n-1)$ matrices given by 
\begin{equation} \label{def:Rkl}
{\cal R}_{kl}=m_{k+1}\vr_{k+1}\delta_{kl}-\frac{m_{k+1}m_{l+1}\vr_{k+1}\vr_{l+1}}{\Sigma_{\vr}},
\end{equation}
\begin{equation} \label{lag:5b}
\cB_{kl}=\frac{\vr_{k+1}\vr_{l+1} D_{k+1,l+1}}{p}.
\end{equation}
for $k,l=1,\ldots,n-1$. 
Moreover, the matrix $\cR$ is uniformly coercive in $(x,t)$ and the same property holds for $\cB$ provided 
that either: \\
{\emph{Condition 1:}} The matrix $C$ is of the form \eqref{Cform}\\
or \\
{\emph{Condition 2:}} $\Omega$ is bounded and \eqref{prop_D} is satisfied for $x\in \Ov{\Omega}$, $t\in[0,T]$.\\

\end{prop}

The local well-posedness  of system \eqref{sys:normal},\eqref{bc:normal} in the maximal $L_p-L_q$ regularity setting is provided by our main result below.
\begin{thm}\label{thm:main2}
Assume that 
\begin{itemize}
\item $2 < p < \infty$, $3 < q < \infty$, $2/p + 3/q < 1$  and $L > 0$;
\item $\Omega$ is a uniform $C^3$ domain in
$\R^3$;
\item there exists a constant $C>0$ such that
\begin{equation} \label{nablaD}
 \forall \, k,l \in 1,\ldots,n \quad \|\nabla D_{kl}(t,\cdot)\|_{L_q(\Omega)}\leq C \sum_{j=1}^3\|\nabla \vr_j(t,\cdot)\|_{L_q(\Omega)} \quad \textrm{a.e. in} \; (0,T);
\end{equation}
\item there exist
positive numbers $a_1$ and $a_2$ for which
\begin{equation}\label{initial:0}
a_1 \leq \vr_{k}^0(x) \leq a_2 \quad
\forall x \in \overline{\Omega}, \; k \in 1, \ldots, n.
\end{equation}
\end{itemize}
Let $\vr_{k}^{0}(x), k=1,\ldots n$, and 
$\bu^0(x)$ be initial data for Eq. \eqref{1.1} 
and let 
$$
(\vr^0(x),h_1^0(x),\ldots,h_{n-1}^0(x)) = \Psi(\vr_1^0(x),\ldots \vr_n^0(x)).
$$
Then, there exists a time $T>0$ depending on
$a_1$, $a_2$ and $L$ such that if 
the initial data satisfy the condition:
\begin{equation}\label{initial:1}
\|\nabla(\vr_1^0,\ldots ,\vr_n^0)\|_{L_q(\Omega)}
+ \|\bu^0\|_{B^{2(1-1/p)}_{q,p}(\Omega)^3} 
+ \|h_1^0,\ldots,h_{n-1}^0\|_{B^{2(1-1/p)}_{q,p}(\Omega)^{n-1}}
\leq L
\end{equation}
and the compatibility condition:
\begin{equation}\label{initial:2}
\bu^0|_\Gamma=0, \quad \nabla h^0_{k} \cdot \bn|_\Gamma = 0, \quad k=1,\ldots,n-1,
\end{equation}
then problem \eqref{sys:normal} with boundary conditions \eqref{bc:normal}   and initial conditions \eqref{ic:normal} admits a unique solution 
$(\vr, \bu, h_1,\ldots,h_{n-1})$ with
\begin{gather*}
\vr - \vr^0 \in H^1_p((0, T), H^1_q(\Omega)),
\quad \bu \in H^1_p((0, T), L_q(\Omega)^3) \cap L_p((0, T), H^2_q(\Omega)^3),\\
h_1,\ldots,h_{n-1} \in H^1_p((0, T), L_q(\Omega)) \cap L_p((0, T), H^2_q(\Omega))
\end{gather*}
possessing the estimates:
\begin{gather*}
\|\vr-\vr^0\|_{H^1_p((0, T), H^1_q(\Omega))}
+ \|\pd_t(\bu, h_1,\ldots,h_{n-1})\|_{\lpqtn}
+ \|(\bu, h_1,\ldots,h_{n-1})\|_{L_p((0, T), H^2_q(\Omega)^{n+2})}
\leq CL, \\ a_1 \leq \vr(x,t) \leq na_2+a_1
\quad\text{for $(x, t) \in \Omega\times(0, T)$}, \quad 
 \int^T_0\|\nabla\bu(\cdot, s)\|_{L_\infty(\Omega)}
\leq \delta. 
\end{gather*}
Here, $C$ is some constant independent of $L$, and $\delta$ is a small positive parameter.
\end{thm}
Let us state some remarks concerning our main result.
\begin{rmk}
Notice that due to \eqref{def_D} the requirement \eqref{nablaD} is satisfied for the special form \eqref{Cform} provided $C_1 \leq |\vr_k| \leq C_2$ for some positive constants $C_1<C_2$. 
\end{rmk}
\begin{rmk}
The parameter $\delta$ above remains small for large times. This is especially important for the existence of global-in-time solutions, not included in the present study.
\end{rmk}

\begin{rmk}
Due to conditions \eqref{coerc:B} we can apply the inverse of $\cB$ to the boundary conditions \eqref{bc:normal} which leads to equivalent formulation of the boundary condition in the standard form
\begin{equation} \label{bc:normal1}
\bu=0, \quad \nabla h_{k} \cdot \bn = 0, \quad k=1,\ldots,n-1,\quad\mbox{on}\ (0,T)\times\partial\Omega .
\end{equation}
\end{rmk}
\begin{rmk}
The condition $\frac{2}{p}+\frac{3}{q}<1$ deserves a more detailed comment. First of all, it is stronger than condition $\frac{2}{p}+\frac{3}{q} \neq 1$ imposed 
in Theorems \ref{thm:lin1} and \ref{thm:main1}, which gives solvability of associated linear problems. A natural question is whether the condition in Theorem \ref{thm:main2} cannot be strenghtened. The answer is partially positive. One could relax this condition allowing $\frac{2}{p}+\frac{3}{q}>1$ with additional constraints on $p,q$ following \cite{SSZ}. However, this would be at a price of numerous additional technicalities that we omit here for brevity. 
\end{rmk}

A keynote requirement necessary to prove our main result is the coercivity of matrices ${\cal R}$ and $\cB$. The details are given in the Appendix, however it is worth to mention in this place that we need to know that fractional densities are bounded from below by a positive constant. Note that the statement of Theorem \ref{thm:main2} provides us only  with bounded  functions $h_i$ given by \eqref{def:psi}. Let us therefore check that these conditions are in fact equivalent. The implication in one direction follows immediately from \eqref{def:psi}, for the other one we have:
\begin{lem}\label{lem:1}
Let $h_i$ given by \eqref{def:psi} be bounded and let
\begin{equation} \label{vrpos:0}
\vr \geq C>0.
\end{equation}
Then
\begin{equation}\label{rhoidown}
\vr_i \geq C>0, \quad i=1,\ldots,n.
\end{equation}
\end{lem}
\emph{Proof.} Assume $\exists  i \in\{1,\ldots, n-1\}$ and $(x_0,t_0)$ s.t. 
$$
\lim_{(x,t)\longrightarrow (x_0,t_0)}\vr_{i+1}(x,t)=0.
$$
Then 
\begin{equation} \label{vrpos:1}
\lim_{(x,t)\longrightarrow (x_0,t_0)}\vr_1(x,t)=0 
\end{equation}
since otherwise $h_i(x,t)$ would be unbounded from below. This in turn implies 
that 
\begin{equation} \label{vrpos:2}
\lim_{(x,t)\longrightarrow (x_0,t_0)}\vr_{k+1}(x,t)=
 0 \quad \forall\; 1 \leq k \leq n-1
\end{equation}
since otherwise corresponding $h_{k}$ would be unbounded from above. This means that $\sum_{k=1}^3\vr_k(x,t)=0$ which  contradicts \eqref{vrpos:0}.

\qed

Let us finish this section with presenting the outline of the rest of the paper. 
In Section \ref{S:Lag} we  rewrite the problem in Lagrangian coordinates ; this step is necessary to apply 
the maximal $L_p-L_q$ regularity theory. In Section \ref{S:lin} we linearize the problem around the initial 
condition.  Section \ref{S:Nonl}  is dedicated to nonlinear estimates which are used to close the fixed point 
argument and prove Theorem \ref{thm:main2} using the existence result for linearized system from  
Theorem \ref{thm:main1}, the proof of which can be found in \cite{PSZ2}.   

\section{Lagrangian coordinates}\label{S:Lag}
We begin the proof of Theorem \ref{thm:main2} by transforming the symmetrized system \eqref{sys:normal} to the Lagrangian coordinates $x = \Phi(y,t)$ related to the vector field $\bv$:  
\begin{equation}\label{lag:1}
x = y + \int^t_0\bv(y, s)\,ds.
\end{equation}
Then for any differentiable function $f$ we have 
\begin{equation} \label{dt_lag}
\de_t f(\Phi(t,y),t)=\de_t f+\vu \cdot \nabla_x f.
\end{equation}
Since
\begin{equation}\label{lag:2}
\frac{\pd x_i}{\pd y_j} = \delta_{ij} 
+ \int^t_0\frac{\pd v_i}{\pd y_j}(y, s)\,ds,
\end{equation}
assuming that 
\begin{equation}\label{assump:1}
\sup_{t \in (0,T)}\int^t_0\|\nabla\bv(\cdot, s)\|_{L_\infty(\Omega)}\,ds
\leq \delta
\end{equation}
for sufficiently small positive constant $\delta$, the matrix
$\pd x/\pd y = (\pd x_i/\pd y_j)$ has the inverse 
\begin{equation}\label{lag:3}
\Bigr(\frac{\pd x_i}{\pd y_j}\Bigr)^{-1} = \bI + \bV^0(\bk_{\bv}), \quad \bk_{\bv} = \int^t_0\nabla\bv(y, s)\,ds.
\end{equation}
Here , $\bI\,$ is the $3\times 3$
identity matrix, and $\bV^0(\bk)$ is the $3\times 3$ matrix of 
smooth functions with $\bV^0(0) = 0$. We have
\begin{equation}\label{lag:4}
\nabla_x = (\bI + \bV^0(\kv))\nabla_y, 
\quad \frac{\pd}{\pd x_i} = \sum_{j=1}^3 (\delta_{ij} + V^0_{ij}(\kv))
\frac{\pd}{\pd y_j}.
\end{equation}
Moreover (see for instance \cite{St1}), the map $\Phi(y, t)$ is bijection from $\Omega$ onto
$\Omega$.

We define our unknown functions in Lagrangian coordinates:
\begin{equation}\label{lag:5}
\bv(y, t) = \bu(x, t),
\quad \eta(y, t) = \vr(x, t), \quad
\vt_i(y, t) = h_i(x, t), \; i=1,\ldots,n-1, 
\end{equation}
and we denote 
$$\vec{\vt}:=(\vt_1,\ldots,\vt_{n-1})^\top.$$ 
We now show that $U=(\bv,\eta,\vec\vt)$ satisfies the system  
\eq{\label{lag:sys}
&\de_t\eta + \eta\dv\bv = R_1(U) 
\\
&\eta \de_t \bv - \mu \Delta \bv - (\mu+\nu)\nabla\dv\bv +\frac{\eta}{\Sigma_\vr}\nabla\eta
+\sum_{l=1}^{n-1}\lr{\vr_{l+1}-\frac{m_{l+1}\vr_{l+1}\vr}{\Sigma_\vr}}\nabla \vt_l = \vc{R}_2(U)\\
&\sum_{l=1}^{n-1} {\cal R}_{kl}\de_t \vt_l + \lr{\vr_{k+1}-\frac{m_{k+1}\vr_{k+1}\vr}{\Sigma_\vr}}\Div \bv-\dv\lr{\sum_{l=1}^{n-1} \cB_{kl}\nabla\vt_l}=R^k_3(U), \quad k=1,\ldots,n-1
}
supplemented with the  boundary conditions  
\begin{equation} \label{lag:bc}
\bv|_{\de \Omega}=0, \quad \nabla \vt_k\cdot\bn|_{\de \Omega}=R^k_4(U), \quad k=1,\ldots,n-1
\end{equation}
where 
\begin{equation} \label{def:vrk}
(\vr_1,\ldots,\vr_n)=(\vr_1,\ldots,\vr_n)(\eta,\vec{\vt})=\Psi^{-1}(\eta,\vec{\vt}).
\end{equation}
\begin{rmk}
In the remainder of the paper we write simply $\vr_k$ keeping in mind that we have 
the dependence \eqref{def:vrk} since we work in Lagrangian coordinates. 
\end{rmk}
We now derive the  precise form of terms on the right hand side of \eqref{lag:sys},\eqref{lag:bc}.
First of all we have 
\begin{equation}\label{lag:div}
\dv_x = \dv_y + \sum_{i,j=1}^3V^0_{ij}(\kv)\frac{\de v_i}{\de y_j},
\end{equation} 
therefore we easily obtain \eqref{lag:sys}$_1$ with 
\begin{equation}\label{lag:6}
R_1(U) = -\eta\sum_{i,j=1}^3 V^0_{ij}(\kv)\frac{\pd v_i}{\pd y_j}.
\end{equation} 
Now we need to transform second order operators.   
By \eqref{lag:4}, we have
$$
\Delta_x \bu = \sum_{k=1}^3\frac{\pd}{\pd x_k}\lr{\frac{\pd \bu}{\pd x_k}}
= \sum_{k,l,m=1}^3\lr{\delta_{kl} + V^0_{kl}(\kv)}
\frac{\pd}{\pd y_l}
\lr{\lr{\delta_{km} + V^0_{km}(\kv)}\frac{\pd \bv}{\pd y_m}}.
$$
Therefore 
$$\Delta_x \bu = \Delta_y \bv + A_{2\Delta}(\kv)\nabla^2_y\bv
+ A_{1\Delta}(\kv)\nabla_y\bv
$$
with
\eq{
A_{2\Delta}(\kv)\nabla^2_y\bv &= 2\sum_{l,m=1}^3V^0_{kl}(\kv)
\frac{\pd^2\bv}{\pd y_l\pd y_m}
+ \sum_{k,l, m=1}^3
V^0_{kl}(\kv)V^0_{km}(\kv)
\frac{\pd^2\bv}{\pd y_l \pd y_m}, \label{a2delta} }
\eq{
A_{1\Delta}(\kv)\nabla_y\bv  = &\sum_{l, m=1}^3(\nabla_\kv V^0_{l m})(\kv)
\int^t_0(\pd_l\nabla_y\bv)\,ds \frac{\pd \bv}{\pd y_m}\\
&+ \sum_{k,l, m=1}^3
V^0_{kl}(\kv) (\nabla_\kv V^0_{km})(\kv)
\int^t_0\pd_l\nabla_y\bv\,ds\frac{\pd \bv}{\pd y_m}, \label{a1delta}
}
where $(\nabla_\kv V^0_{km})(\kv)$ denotes $ \lr{V^0_{km}}'(\kv)$.

Similarly for $i\in\{1,\ldots,N\}$ we have
$$\frac{\pd}{\pd x_j}\dv_x\bu 
= \sum_{k=1}^3(\delta_{jk} + V^0_{jk}(\kv))\frac{\pd}{\pd y_k}
\lr{\dv_y\bv + \sum_{l, m=1}^3 V^0_{l m}(\kv)\frac{\pd v_l}{\pd y_m}},
$$
so we obtain
$$\frac{\pd}{\pd x_j}\dv_x\bu
= \frac{\pd}{\pd y_j}\dv_y\bv + A_{2\dv, j}(\kv)\nabla^2_y\bv
+ A_{1\dv, j}(\kv)\nabla_y\bv,
$$
where
\eq{
A_{2\dv,j}(\kv)\nabla^2_y\bv
& = \sum_{l, m=1}^3V^0_{l m}(\kv)\frac{\pd^2 v_l}{\pd y_m \de y_j}
+ \sum_{k=1}^3 V^0_{jk}(\kv)\frac{\pd}{\pd y_k}\dv_y\bv
+ \sum_{k, l=1}^3V^0_{jk}(\kv)V^0_{l m}(\kv)
\frac{\pd^2v_l}{\pd y_k \pd y_m}, \label{a2div}
}
\eq{
A_{1\dv, j}(\kv)\nabla_y\bv
 =& \sum_{l, m=1}^3(\nabla_{\kv} V^0_{l m})(\kv)
\int^t_0\pd_j\nabla_y\bv\,ds\frac{\pd v_l}{\pd y_m} \\
&+ \sum_{k,l, m=1}^3V^0_{jk}(\kv)(\nabla_{\kv} V^0_{l m})(\kv)
\int^t_0\pd_k\nabla_y\bv\,ds\frac{\pd v_l}{\pd y_m}. \label{a1div}
}
Therefore, transforming also $\nabla_x \vr$ and $\nabla_x h_l$ we obtain \eqref{lag:sys}$_2$   
with
\eq{ \label{lag:7} 
\vc{R}_2(U) =& \mu A_{2\Delta}(\kv)\nabla^2_y\bv 
+ \mu A_{1\Delta}(\kv)\nabla_y\bv
+ \nu A_{2\dv}(\kv)\nabla^2_y\bv + \nu A_{1\dv}(\kv)\nabla_y\bv \\
&+ \frac{\eta}{\Sigma_\vr}\bV^0(\kv)\nabla_y\eta
+ \bV^0(\kv)\sum_{l=2}^3\lr{\vr_l-\frac{m_l\vr_l\vr}{\Sigma_\vr}}\nabla_y \vt_{l-1},
}
where  $A_{i\dv}\nabla^i_y\bv,\  i=1,2$ are vectors with coordinates $A_{i\dv,j}\nabla^i_y\bv,\  j=1,\ldots,N$. 

Finally we transform the species balance equations. 
We have 
\begin{equation*}\begin{split}
\dv_x(\cB_{kl}\nabla_x h_l)
&=\cB_{kl}(\Delta_y\vt_l+A_{2\Delta}(\kv)\nabla^2_y\vt_l+A_{1\Delta}(\kv)\nabla_y \vt_l)\\
&\quad+\lr{\nabla_y \cB_{kl}+\bV^0(\kv)\nabla_y \cB_{kl}}\lr{\nabla_y \vt_l+\bV^0(\kv)\nabla_y \vt_l}\\
&=\dv_y(\cB_{kl}\nabla_y\vt_l)+R^{kl}_3(U),
\end{split}\end{equation*}
where 
\eq{\label{lag:8} 
R^{kl}_3(U)=&\cB_{kl}(A_{2\Delta}(\kv)\nabla^2_y\vt_l+A_{1\Delta}(\kv)\nabla_y \vt_l)\\
&+\bV^0(\kv)\nabla_y \cB_{kl}(\nabla_y\vt_l+\bV^0(\kv)\nabla_y\vt_l)+(\nabla_y \cB_{kl})\bV^0(\kv)\nabla_y\vt_l.
}
Therefore, transforming also $\div \bu$, we obtain \eqref{lag:sys}$_3$ with  
\begin{equation}\label{lag:10}  
R^k_{3}(U)=\sum_{l=1}^{n-1} R_{3}^{kl}(U)-\lr{\vr_{k+1}-\frac{m_{k+1}\vr_{k+1}\vr}{\Sigma_\vr}}\sum_{j,m=1}^3V^0_{jm}(\kv)\frac{\de v_j}{\de y_m}.  
\end{equation}  
It remains to transform the boundary conditions. For this purpose notice that 
$$\bn(x) = \bn\lr{y + \int^t_0\bv(y, s)\,ds}
= \bn(y) + \int^1_0(\nabla\bn)
\lr{y + \tau\int^t_0\bv(y, s)\,ds}\,d\tau
\int^t_0\bv(y, s)\,ds,
$$
and therefore we obtain \eqref{lag:bc} with  
\eq{\label{lag:9}
R_4^k(U)=&\bn\lr{y + \int^t_0\bv(y, s)\,ds}\cdot (\bV^0(\kv)\nabla_y \vartheta_k)\\
&+ \left\{\int^1_0(\nabla\bn)
\lr{y + \tau\int^t_0\bv(y, s)\,ds}\,d\tau
\int^t_0\bv(y, s)\,ds\right\} \cdot\nabla_y\vartheta_k.
}

\section{Linearization}\label{S:lin}
\subsection{Formulation of linearized system}
We now  linearize the system in the Lagrangian coordinates \eqref{lag:sys} around the initial conditions.
For this purpose we introduce small perturbations 
\begin{equation} \label{lin1:1}   
\eta=\sigma+\vr^0,\quad \vr_l=\sigma_l+\vr^0_l,
\end{equation}  
following the convention introduced in the previous section that $\vr_l$ are the functions in the Lagrangian coordinates.

\noindent Let us denote 
$$
\Sigma_\vr^0=\sum_{k=1}^n m_k \vr^0_k, \quad p^0=\sum_{k=1}^n \frac{\vr_k^0}{m_k},
$$
and
\begin{equation} \label{lin1:1b}
h^0_k=\frac{1}{m_k}\log \vr^0_{k+1} - \frac{1}{m_1}\log \vr^0_1, \quad k=1,\ldots,n-1 .
\end{equation}
Observe that due to \eqref{initial:0} we have
\begin{equation}
n a_1 \leq \vr^0 \leq n a_2,
\end{equation}
as well as 
$$  
|h^0_k|\leq \frac{1}{m_{k+1}}|\log a_2|+\frac{1}{m_1}|\log a_1|.  
$$

The linearization of the continuity equation is straighforward, while for the momentum equation we have 
$$
\frac{\eta}{\Sigma_{\vr}}\nabla \eta = \frac{\vr^0}{\Sigma_\vr^0} \nabla \sigma 
+ \vr^0 \nabla \sigma \left( \frac{1}{\Sigma_\vr}-\frac{1}{\Sigma_\vr^0} \right) + \frac{\eta}{\Sigma_\vr}\nabla\vr^0
+ \frac{\sigma}{\Sigma_\vr}\nabla\sigma 
$$  
and 
\begin{equation} \label{lin1:2}
\frac{m_l\vr_l\vr}{\Sigma_\vr}=
\frac{m_l\vr^0_l\vr^0}{\Sigma_\vr^0}+m_l\vr^0\vr^0_l\left( \frac{1}{\Sigma_\vr}-\frac{1}{\Sigma_\vr^0} \right)
+\frac{m_l}{\Sigma_\vr}(\vr_l^0\sigma+\vr^0\sigma_l). 
\end{equation}
Similarly we linearize the $\cR_{kl}$ in the species equations
while for the reduced diffusion matrix we use
\begin{equation} \label{lin1:3}
\cB_{k-1,l-1}=\frac{\vr_k \vr_l D_{kl}}{p}=\frac{\vr^0_k\vr^0_lD_{kl}^0}{p^0}
+\frac{\vr_k\vr_lD_{kl}-\vr^0_k\vr^0_lD_{kl}^0}{p}+\vr^0_k\vr^0_lD_{kl}^0\left(\frac{1}{p}-\frac{1}{p^0}\right).
\end{equation}
Therefore we obtain the following linearized system   
\begin{align} \label{lin1:sys}
&\pt\sigma + \vr^0 \dv \bv = f_1(U)\\
&\vr^0 \de_t \bv - \mu\Delta\bv-(\mu+\nu)\nabla\dv \bv + \gamma_1 \nabla \sigma 
+ \sum_{l=1}^{n-1} \gamma_2^{l} \nabla\vt_{l}=\bbf_2(U)\\
&\sum_{l=1}^{n-1} {\cal R}_{kl}^0\de_t \vt_l + \gamma_2^k\Div \bv
-\dv\lr{\sum_{l=1}^{n-1} \cB_{kl}^0\nabla\vt_l}=f^k_3(U), \quad k=1,\ldots,n-1
\end{align}
 in $\Omega\times(0, T)$, supplied with the boundary conditions 
\begin{equation} \label{lin1:bc}
\bv|_{\de \Omega}=0, \quad \nabla \vt_k \cdot \bn|_{\de \Omega}=f^k_4(U), \quad k=1,\ldots,n-1
\end{equation}
and initial conditions 
\begin{equation} \label{lin1:ic}
(\sigma,\bv,\vec\vt)|_{t=0}=(0,\bu^0, \vec h^0),
\end{equation}
where we denote
$$\vec h^0=(h_1^0,\ldots,h_{n-1}^0),$$
$$D_{kl}^0=D_{kl}(\vr^0), \quad \cR _{kl}^0=m_{k+1}\vr^0_{k+1}\delta_{kl}-\frac{m_{k+1}m_{l+1}\vr^0_{k+1}\vr^0_{l+1}}{\Sigma_{\vr}^0},  \quad \cB_{kl}^0=\frac{\vr^0_{l+1}\vr^0_{k+1}D_{k+1,l+1}^0}{p_0},$$
$$\gamma_1=\frac{\vr^0}{\Sigma_\vr^0}, \quad   \gamma_2^k=\vr^0_{k+1}-\frac{m_{k+1}\vr^0_{k+1}\vr^0}{\Sigma_\vr^0}, $$
and the right hand side is given by 
\begin{equation} \label{lin1:5a} 
f_1(U)=R_1(U)-\sigma \dv \bv,
\end{equation}  
\begin{equation} \label{lin1:5b}
\begin{split}
\bbf_2(U)=&R_2(U)-\sigma\de_t \bv - \vr^0\nabla\eta \left( \frac{1}{\Sigma_\vr}-\frac{1}{\Sigma_\vr^0} \right) - \frac{\vr^0}{\Sigma_\vr}\nabla\vr^0
- \frac{\sigma}{\Sigma_\vr}\nabla\eta\\
&+\sum_{l=1}^{n-1}\lr{-\sigma_{l+1}+m_{l+1}\vr_{l+1}^0\vr^0\left( \frac{1}{\Sigma_\vr}-\frac{1}{\Sigma_\vr^0} \right)+\frac{m_{l+1}}{\Sigma_\vr}(\vr_{l+1}\sigma+\vr^0\sigma_{l+1})}\nabla\vt_{l},
\end{split}
\end{equation}
\eq{ \label{lin1:5c}
f^k_3(U)
&=R^k_3(U) + \vr^0_{k+1}\dv \bv + \left[ m_{k+1}\vr^0\vr^0_{k+1}\lr{\frac{1}{\Sigma}-\frac{1}{\Sigma^0}}+\frac{m_{k+1}}{\Sigma_\vr}(\vr_{k+1}^0\sigma+\vr^0\sigma_{k+1})\right]\dv \bv\\
&+\sum_{l=1}^{n-1}\left( - \delta_{kl}m_{k+1}\sigma_{k+1}
+m_{k+1} m_{l+1}\left[\vr^0_{k+1}\vr^0_{l+1}\left(\frac{1}{\Sigma}-\frac{1}{\Sigma_\vr^0}\right)+
\frac{\vr_{k+1}^0\sigma_{l+1}+\vr^0_{l+1}\sigma_{k+1}}{\Sigma_\vr}
\right]\right)\de_t\vt_{l}\\
&+\dv \left( \sum_{l=1}^{n-1} \left[\frac{\vr_{k+1}\vr_{l+1}D_{k+1,l+1}-\vr^0_{k+1}\vr^0_{l+1}D_{k+1,l+1}^0}{p}+\vr^0_{k+1}\vr^0_{l+1}D_{k+1,l+1}^0\lr{\frac{1}{p}-\frac{1}{p_0}}\right]\nabla\vt_{l} \right),
}
\begin{equation} \label{lin1:5d}
f^k_4(U)=R^k_4(U).
\end{equation}

\subsection{Solvability of the complete linear system} 
\subsubsection{Auxiliary results}
To prove existence of local-in-time strong solutions to system \eqref{lin1:sys} with fixed and given  $f_1, \bbf_2, f_3^k$, and $f_4^k$ we will use some auxiliary results for two subsystems. First let us recall a relevant existence result for the  fluid part (for the proof see \cite{PSZ}, Theorem 5.1 ):
\begin{thm} \label{thm:lin1}
Assume $1 < p, q < \infty$  
$2/p + 1/q \not =1$, $T>0$ and
$\Omega$ is a uniformly $C^2$ domain in $\BR^N$ $(N \geq 2)$.
Assume moreover that $\vr^0\in H^1_q(\Omega)$, 
$\bu_0 \in B^{2(1-1/p)}_{q,p}(\Omega)^N$, 
$\tilde f_1 \in L_p(\BR, H^1_q(\Omega)^N) $ and $\tilde \bbf_2 \in L_p((0, T), L_q(\Omega)^N)$. 
Then the problem
\begin{equation}
\left\{
\begin{aligned}
\pt\sigma+ \vr^0 \dv \bv &= \tilde f_1 &\quad&\text{in $\Omega\times(0, T)$}, \\
\vr^0 \de_t \bv - \mu\Delta\bv-(\mu+\nu)\nabla\dv \bv + \gamma_1 \nabla \sigma &= \tilde\bbf_2&\quad&\text{in $\Omega\times(0, T)$}, \\
\bv|_{\de \Omega}&=0&\quad&\text{on $\Gamma \times (0, T)$}, \\
(\sigma,\bv)|_{t=0}&=(0,\bu^0)&\quad&\text{in $\Omega$},
\end{aligned}\right.
\end{equation}
admits a solution $(\sigma,\bv)$ such that
\begin{align}
&\|\bv\|_{L_p((0, T), H^2_q(\Omega))}
+ \|\pd_t\bv\|_{L_p((0, T), L_q(\Omega))}
+ \| \sigma \|_{H^1_p(0,T;H^1_q(\Omega))} \\
&\quad
\leq Ce^{cT}\lr{\|\vr^0\|_{H^1_q(\Omega)}+\|\bu^0\|_{B^{2(1-1/p)}_{q,p}(\Omega)}
+ \|\tilde f_1\|_{L_p((0, T), H^1_q(\Omega))}
+ \|\tilde\bbf_2\|_{L_p((0, T), L_q(\Omega))}}.
\end{align}
\end{thm}
For the species subsystem we recall the following theorem which gives solvability in a maximal $L_p-L_q$ 
regime of a linear problem, its proof can be found in our previous work \cite{PSZ2}.
For general $m$ species we consider $k\in\{1,\ldots,m\}$ and the following set of equations
\begin{equation}\label{1.1?}\left\{
\begin{aligned}
\sum_{\ell=1}^m \cR_{k\ell}\pd_t \vt_\ell
-\dv\lr{\sum_{\ell=1}^m \cB_{k\ell}\nabla \vt_\ell} & = \tilde f_3^k
&\quad&\text{in $\Omega\times(0, T)$}, \\
\sum_{\ell=1}^m \cB_{k\ell}\nabla \vt_\ell \cdot \bn & = \tilde f_4^k
&\quad&\text{on $\Gamma \times (0, T)$}, \\
\vt_k|_{t=0} & = h_{k}^0
&\quad&\text{in $\Omega$},
\end{aligned}\right.
\end{equation}
where 
$\cB= \cB(x)$ and $\cR=\cR(x)$ are $m\times m$ matrices whose
$(k, \ell)^{\rm th}$ components are $\cB_{k\ell}(x)$ and $\cR_{k\ell}(x)$, respectively.
\begin{thm} \label{thm:main1}
Assume that 
\begin{itemize}
\item
there exists a number  $M_0$ for which
\begin{equation}\label{1.2}\begin{aligned}
&|\cB_{k\ell}(x)|, |\cR_{k\ell}(x)| \leq M_0, 
\quad \text{for any $x \in \Omega$}, \\
 &|\cB_{k\ell}(x) - \cB_{k\ell}(y)|\leq M_0|x-y|^\sigma,
\quad
|\cR_{k\ell}(x) - \cR_{k\ell}(y)|\leq M_0|x-y|^\sigma 
\quad\text{for any $x, y \in \Omega$},\\
&\|\nabla(\cB_{k\ell}, \cR_{k\ell})\|_{L_r(\Omega)} \leq M_0.
\end{aligned}
\end{equation}
\item
the matrices
$\cB$ and $\cR$ are positive and symmetric and 
that there exist constants $m_1,m_2 > 0$ for which
\begin{equation}\label{1.3}
\langle \cB(x)\xi, \overline{\xi} \rangle \geq m_1|\xi|^2,
\quad 
\langle \cR(x)\xi, \overline{\xi} \rangle \geq m_2|\xi|^2
\end{equation}
for any complex $m$-vector $\xi$ and $x \in \Omega$.

\item
$1 < p, q < \infty$ and $T > 0$, 
$2/p + 1/q \not =1$  and
$\Omega$ is a uniformly $C^2$ domain in $\BR^3$ $(N \geq 2)$.

\item for all $k=1,\ldots,m$, 
$h_k^0 \in B^{2(1-1/p)}_{q,p}(\Omega)$, 
$\tilde f_3^k \in L_p((0, T), L_q(\Omega))$ and 
$\tilde f_4^k \in L_p(\BR, H^1_q(\Omega)) 
\cap
H^{1/2}_p(\BR, L_q(\Omega))$ are given functions satisfying
the compatibility conditions:
\begin{equation}\label{1.4}\sum_{\ell=1}^m B_{k\ell}\nabla h_{\ell}^0
\cdot\bn = \tilde f_4^k(\cdot, 0)
\quad\text{on $\Gamma$}
\end{equation}
provided $2/p + 1/q < 1$.  
\end{itemize}

\noindent
Then, problem \eqref{1.1?}
admits a unique solution $\vec\vt = (\vt_1, \ldots, \vt_m)^\top$
with
\begin{equation}\label{1.5}
\vec\vt \in L_p((0, T), H^2_q(\Omega)^m)
\cap H^1_p((0, T), L_q(\Omega)^m)
\end{equation}
possessing the estimate:
\begin{equation}\label{1.6}
\begin{aligned}
&\|\vec\vt\|_{L_p((0, T), H^2_q(\Omega))}
+ \|\pd_t\vec\vt\|_{L_p((0, T), L_q(\Omega))}\\
&\quad
\leq Ce^{cT}(\|\vec h^0\|_{B^{2(1-1/p)}_{q,p}(\Omega)}
+ \|\vec{\tilde f}_3\|_{L_p((0, T), L_q(\Omega))}
+ \|\vec{\tilde f}_4\|_{L_p((0, T), H^1_q(\Omega))}
+ \|\vec{\tilde f}_4\|_{H^{1/2}_p(\BR, L_q(\Omega))})
\end{aligned}
\end{equation}
for some constants $C$ and $c$.
\end{thm}

\subsubsection{Fixed point argument}

With Theorems \ref{thm:lin1} and \ref{thm:main1} it is easy to show solvability with appropriate estimates of complete linear system corresponding to \eqref{lin1:sys}-\eqref{lin1:bc}:
\begin{equation}\label{lin2:sys}
\left\{
\begin{aligned}
&\pt\sigma + \vr^0\dv \bv = f_1 \\
&\vr^0 \de_t \bv - \mu\Delta\bv-(\mu+\nu)\nabla\dv \bv + \gamma_1 \nabla \sigma +\sum_{l=1}^{n-1}\gamma^l_2 \nabla \vt_l = \bbf_2\\
&\sum_{l=1}^{n-1} {\cal R}_{kl}^0\de_t \vt_l + \gamma_2^k\Div \bv
-\dv\lr{\sum_{l=1}^{n-1} \cB_{kl}^0\nabla\vt_l}=f^k_3, \quad k=1,\ldots,n-1
\end{aligned}\right.
\end{equation}
with given $\gamma_1,\{\gamma^l_2\}_{l=1,\ldots,n-1}$  and the boundary conditions 
\begin{equation} \label{lin2:bc}
\bv|_{\de \Omega}=0, \quad \sum_{l=1}^{n-1} \cB_{kl}^0 \nabla \vt_l\cdot \bn|_{\de \Omega}=f_4^k, \quad k=1,\ldots,n-1.
\end{equation}
and initial conditions \eqref{lin1:ic}. 

We have the following result.
\begin{thm} \label{thm:lin2}
Assume $\cB^0$,$\cR^0$, $\Omega$ and $p,q$ satisfy the assumptions of Theorem \ref{thm:main1} with $m=n-1$. Assume moreover 
$\bu^0,\vec h^0 \in B^{2(1-1/p)}_{q,p}(\Omega)$,
$\vr^0 \in H^1_q(\Omega)$,
$f_1 \in L_p((0, T), H^1_q(\Omega))$, 
$(\bbf_2,\vec f_3) \in L_p((0, T), L_q(\Omega)^{n+2})$,
$\vec f_4 \in L_p(\BR, H^1_q(\Omega)^{n-1}) 
\cap
H^{1/2}_p(\R, L_q(\Omega)^{n-1})$. 
Then for any $M>0$, if 
\begin{align}
&\|\bu^0,\vec h^0\|_{B^{2(1-1/p)}_{q,p}(\Omega)}
+\|\vr^0\|_{H^1_q(\Omega)}
+\|f_1\|_{L_p((0, T), H^1_q(\Omega))}\\
&+\|(\bbf_2,\vec f_3)\|_{L_p((0, T), L_q(\Omega)^{n+2})}
+\|\vec f_4\|_{L_p(\BR, H^1_q(\Omega)^{n-1})} 
+\|\vec f_4\|_{H^{1/2}_p(\R, L_q(\Omega)^{n-1})} \leq M,
\end{align}
then there exists $T>0$ such that system
 \eqref{lin2:sys}-\eqref{lin2:bc} admits a solution $(\sigma,\bv,\vec\vt)$ on $(0,T)$ with 
\begin{align} \label{est:lin3}
[\sigma,\bv,\vec\vt]_T \leq
&\|\bu^0,\vec h^0\|_{B^{2(1-1/p)}_{q,p}(\Omega)}
+\|\vr^0\|_{H^1_q(\Omega)}
+\|f_1\|_{L_p((0, T), H^1_q(\Omega))}\\
&+\|(\bbf_2,\vec f_3)\|_{L_p((0, T), L_q(\Omega))}
+\|\vec f_4\|_{L_p(\R, H^1_q(\Omega))} 
+\|\vec f_4\|_{H^{1/2}_p(\R, L_q(\Omega)^3)}
\end{align}
\end{thm}
\emph{Proof.} We use the Banach fixed point argument. For given $\bar \bv \in \CH^1_{T,M}$ denote by $\vec\vt(\bar \bv)$
solution to \eqref{lin2:sys}$_3$ with $\bv=\bar \bv$ and boundary condition \eqref{lin2:bc}$_2$. 
Since 
$
\|\bv\|_{L_\infty(0,T,H^1_\infty(\Omega))}\leq CM,
$ (see estimate \eqref{est:04}) 
therefore by Theorem \ref{thm:main1} such solution exists for arbitrary time $T>0$, it is unique and it satisfies 
\eq{
[\vec\vt(\bar\bv)]_{T,1}\leq  &C(T)\Big(\|\vec h^0\|_{B^{2(1-1/p)}_{q,p}(\Omega)}+
\|\vec f_3\|_{\lpqtn}+E(T)\|\bv\|_{L_p(0,T;H^1_q(\Omega)^3)}\\
&\quad+\|\vec f_4\|_{L_p(\R, H^1_q(\Omega)^{n-1})}
+\|\vec f_4\|_{H^{1/2}_p(\R, L_q(\Omega)^{n-1})}\Big)\\
\leq & C(T,M)\left(1+E(T)\|\bar\bv\|_{L_p(0,T;H^1_q(\Omega)^3)}\right).
}
Therefore for $(\bar \bv,\bar \sigma) \in \CH^1_{T,M} \times \CH^2_{T,M}$ we can define 
$(\bv, \sigma)={\mathcal T}(\bar \bv,\bar \sigma)$ as a unique solution of the first two equations  of system \eqref{lin2:sys} with 
$\vec\vt=\vec\vt(\bar \bv)$ and boundary condition \eqref{lin2:bc}$_1$. By Theorem \ref{thm:lin1} we have 
\eq{
[\sigma]_{T,2}+[\bv]_{T,1} \leq &
 C(T)\Big(\|\vr^0\|_{H^1_q(\Omega)}+\|\bu^0\|_{B^{2(1-1/p)}_{q,p}(\Omega)^3}\\
&\quad+ \|f_1\|_{L_p((0, T), H^1_q(\Omega))}
+ \|\bbf_2\|_{L_p((0, T), L_q(\Omega)^3)}+\|\nabla\vec\vt(\bar\bv)\|_{L_p((0, T), L_q(\Omega)^{n-1})}\Big)\\
\leq & C(T,M)\left(1+E(T)\|\bar\bv\|_{L_p(0,T;H^1_q(\Omega)^3)}\right).
}
Moreover, taking different $\bar\bv_1, \bar\bv_2\in \CH^1_{T,M}$ corresponding to the same initial data $\bu^0$, and then subtracting the  for $\vec\vt(\bar \bv_1)$ and $\vec\vt(\bar \bv_2)$
we get 
$$
[\vec\vt(\bar \bv_1)-\vec\vt(\bar \bv_2)]_{T,1}\leq C(M)E(T)[\bar\bv_1-\bar\bv_2]_{T,1}.
$$
Therefore applying Theorem \ref{thm:lin1} to a difference of two solutions we have 
\eqh{
[{\mathcal T}(\bar \bv_1,\bar \sigma_1)-{\mathcal T}(\bar \bv_2,\bar \sigma_2)]_{T,1;T,2}&\leq 
C(M)E(T) [\bar\bv_1-\bar\bv_2]_{T,1}\\
&\leq C(M)E(T)[(\bar\bv_1-\bar\bv_2,\bar\sigma_1-\bar\sigma_2)]_{T,1;T,2}.
}
Therefore for sufficiently small $T$, $\mathcal T$ is a contraction on a set $\CH^1_{T,M}\times\CH^2_{T,M}$,
and applying the Banach fixed point theorem we complete the proof. 

 \qed
\section{Proof of Theorem \ref{thm:main2}}\label{S:Nonl}

\subsection{Nonlinear estimates}
 
The aim of this section is to prove the following proposition which gives the estimate on the right hand side of linearized system in the regularity required in order to apply Theorem \ref{thm:lin2}. We shall use notation introduced at the beginning of Section 5.2. For brevity in this subsection we will not distinguish between scalar and vector valued functions in notation of functional spaces, except for final estimates.  
\begin{prop} \label{prop:est}
Let $\bar U=(\osigma,\ov,\ovt) \in {\mathcal H}_{T,M} $ for given $T,M>0$, where the initial conditions satisfy the assumptions of Theorem \ref{thm:main2}.  
Let $f_1(U),f_2(U),f_3^k(U)$ and $f^k_4(U)$ be given by \eqref{lin1:5a}-\eqref{lin1:5d}, 
where $R_1(U),R_2(U),R_3^k(U)$ and $R^k_4(U)$ are defined in \eqref{lag:6},\eqref{lag:7},\eqref{lag:8}-\eqref{lag:10} and \eqref{lag:9}, respectively. Then 
\eq{\label{est:nonlin} 
&\|f_1(\bar U)\|_{L_p(0,T;W^1_q(\Omega))} + \|f_2(\bar U)\|_{L_p(0,T,L_q(\Omega)^3)}+\|\vec f_3\|(\bar U)\|_{\lpqtn} \\ 
&+\|\vec f^k(\bar U)\|_{L_p(0,T,H^1_q(\Omega)^{n-1})} + \|\vec f_4(\bar U)\|_{H^{1/2}_p(\R,L_q(\Omega)^{n-1})}  \leq C(M,L) E(T).    
}
\end{prop}
Let us start with recalling some auxiliary results. The first one is due to 
Tanabe (cf. \cite{Tanabe} p.10):
\begin{lem} \label{L:int}
Let $X$ and $Y$ be two Banach spaces such that
$X$ is a dense subset of $Y$ and $X\subset Y$ is continuous.
Then for each $p \in (1, \infty)$  
$$H^1_p((0, \infty), Y) \cap L_p((0, \infty), X) 
\subset C([0, \infty), (X, Y)_{1/p,p})$$
and for every $u\in H^1_p((0, \infty), Y) \cap L_p((0, \infty), X)$ we have
$$\sup_{t \in (0, \infty)}\|u(t)\|_{(X, Y)_{1/p,p}}
\leq (\|u\|_{L_p((0, \infty),X)}^p
+ \|u\|_{H^1_p((0, \infty), Y)}^p)^{1/p}.
$$
\end{lem}
Next two results will be needed to estimate the boundary data. For the first one see [Shibata and Shimizu \cite{SS1}, Lemma 2.7]:
\begin{lem}\label{lem:5.1}
Let $1 < p < \infty$, $3 < q < \infty$ and $0 < T \leq 1$.  Assume that 
$\Omega$ is a uniformly $C^2$ domain. Let 
\begin{align*}
f \in H^1_\infty(\BR, L_q(\Omega)) \cap L_\infty(\BR, H^1_q(\Omega)), 
\quad
g \in L_p(\BR, H^1_q(\Omega)) \cap H^{1/2}_p(\BR, L_q(\Omega)).
\end{align*} 
If we assume that $f \in L_p(\BR, H^1_q(\Omega))$ and 
that $f$ vanishes for $t \notin [0, 2T]$ in addition, then we have
\begin{align*}
&\|fg\|_{L_p(\BR, H^1_q(\Omega))} + \|fg\|_{H^{1/2}_p(\BR, L_q(\Omega))}\\
&\quad\leq C(\|f\|_{L_\infty(\BR, H^1_q(\Omega))}
+T^{(q-3)/(pq)}\|\pd_tf\|_{L_\infty(\BR, L_q(\Omega))}^{(1-3/(2q))}
\|\pd_tf\|_{L_p((\BR, H^1_q(\Omega))}^{3/(2q)})
(\|g\|_{_p(\BR, H^1_q(\Omega))} + \|g\|_{H^{1/2}_p(\BR, L_q(\Omega))}).
\end{align*}
\end{lem}
\begin{rmk} \thetag1~ 
The boundary of $\Omega$ 
was assumed to be bounded in \cite{SS1}. 
However, Lemma \ref{lem:5.1} can be proved using Sobolev's 
inequality and complex interpolation theorem, and so 
employing the same argument as that in the proof of
\cite[Lemma 2.7]{SS1}, 
we can prove Lemma \ref{lem:5.1}. \\
\thetag2~ By Sobolev's inequality, $\|fg\|_{H^1_q(\Omega)}
 \leq C\|f\|_{H^1_q(\Omega)}\|g\|_{L_q(\Omega)}$, and so 
the essential part of Lemma \ref{lem:5.1} is the estimate of
$\|fg\|_{H^{1/2}_p(\BR, L_q(\Omega))}$. 
\end{rmk}
The second result has been shown in Shibata and Shimizu \cite{SS2} for $\Omega=\R^3$ and generalized to a uniform $C^2$ domain in Shibata \cite{S17}: 
\begin{lem}\label{lem:5.2} Let $1 < p, q < \infty$. Assume that 
$\Omega$ is a uniform $C^2$ domain.  Then
$$H^1_p(\BR, L_q(\Omega)) \cap L_p(\BR, H^2_q(\Omega))
\subset H^{1/2}_p(\BR, H^1_q(\Omega)), $$
and 
$$\|\nabla f\|_{H^{1/2}_p(\BR, L_q(\Omega))}
\leq C(\|f\|_{L_p(\BR, H^2_q(\Omega))} 
+ \|\pd_t f\|_{L_p(\BR, L_q(\Omega))}).
$$
\end{lem}
Now we show preliminary estimates for functions from the space $\CH_{T,M}$. 
\begin{lem}
Let $\sigma, \bv, \vt_1 \ldots \vt_{n-1} \in \CH_{T,M}$ and let $\kv,\bV^0(\kv)$ be defined in \eqref{lag:3}. Then 
\begin{align}
&\|\bV^0(\kv),\nabla_{\kv}\bV^0(\kv)\|_{L_\infty(\Omega\times(0,T))} \leq C(M,L)E(T), \label{est:01}\\
&{\rm sup}_{t \in (0,T)} \|\sigma(\cdot,t)\|_{H^1_q(\Omega)}\leq C(M,L)E(T), \label{est:02} \\
&{\rm sup}_{t \in (0,T)}\|{\vec \vt}(\cdot,t)-{\vec h^0}\|_{B^{2(1-1/p)}_{q,p}}+{\rm sup}_{t \in (0,T)}\|\bv(\cdot,t)-\bu_0\|_{B^{2(1-1/p)}_{q,p}}\leq C(M,L),\label{est:03}\\
&\|\bv,{\vec \vt}\|_{L_\infty(0,T,H^1_\infty(\Omega))}\leq C(M), \label{est:04}\\
&\|\vr_k-\vr_k^0\|_{L_\infty(0,T;H^1_q)}\leq C(M,L) \quad \forall k=1,\ldots,n,\label{est:05}\\
&\|\vr_k-\vr_k^0\|_{L_\infty(\Omega \times (0,T))} \leq C(M,L)E(T). \label{est:06}
\end{align}
\end{lem}
\emph{Proof}.  First of all, we have 
\begin{align}
\int^T_0\|\nabla\bv(\cdot, t)\|_{L_\infty(\Omega)}\,dt
&\leq C\int^T_0\|\bv(\cdot, t)\|_{H^2_q(\Omega)}\,dt\nonumber\\
&\leq T^{1/{p'}}
\Bigl({\rm sup}_{t \in (0,T)}\int^T_0\|\bv(\cdot, t)\|_{H^2_q(\Omega)}^p\,dt\Bigr)^{1/p}
\leq MT^{1/p'},
\end{align}
which implies \eqref{est:01}. Next,
$$\|\sigma(\cdot, t)\|_{H^1_q(\Omega)}
\leq \int^t_0\|\de_t\sigma(\cdot, s)\|_{H^1_q(\Omega)}\,ds
\leq T^{1/{p'}}\|\de_t\sigma\|_{L_p((0, T), H^1_q(\Omega))}
\leq C(M)E(T),
$$
and so we have \eqref{est:02}. In order to prove \eqref{est:03} we introduce
extension operator
\begin{equation} \label{def:ext} e_T[f](\cdot, t)
= \begin{cases}
0 \quad &t\in(-\infty,0)\cup (2T,+\infty), \\ 
f(\cdot, t) \quad &t\in(0,T), \\
f(\cdot, 2T-t)\quad & t\in(T,2T).
\end{cases}
\end{equation}
Obviously, $e_T[f](\cdot, t) = f(\cdot, t)$ for $t \in (0, T)$.  If
$f|_{t=0}=0$,  then we have
\begin{equation} \label{ext:2} \pd_te_T[f](\cdot, t)
= \begin{cases}
0 \quad &t\in(-\infty,0)\cup (2T,+\infty),\\
(\pd_tf)(\cdot, t) \quad &t\in(0,T), \\
-(\pd_tf)(\cdot, 2T-t)\quad & t\in(T,2T),
\end{cases}
\end{equation}
understood in a weak sense.
Applying Lemma \ref{L:int} with $X=H^2_q(\Omega), \, Y=L_q(\Omega)$ and using 
\eqref{def:ext} and \eqref{ext:2} we have
\begin{align*}
&\sup_{t \in (0, T)}\|\vt(\cdot, t)-h_0\|_{B^{2(1-1/p)}_{q,p}(\Omega)}
\leq \sup_{t \in (0, \infty)}\|e_T[\vt_k-h_k^0]\|_{B^{2(1-1/p)}_{q,p}(\Omega)}\\&\quad 
= (\|e_T[\vt_k-h_k^0]\|_{L_p((0, \infty), H^2_q(\Omega))}^p
+ \|e_T[\vt_k-h^0_k]\|_{H^1_p((0, \infty), L_q(\Omega))}^p)^{1/p}\\
&\quad \leq C(\|\vt_k-h^0_k\|_{L_p((0, \infty), H^2_q(\Omega))}
+ \|\pd_t\vt_k\|_{L_p((0, T), L_q(\Omega))}) \leq C(M,L),
\end{align*}
and estimating $\|\bv(\cdot, t)-\bu_0\|_{B^{2(1-1/p)}_{q,p}(\Omega)}$ 
in the same way we obtain \eqref{est:03}. Then \eqref{est:04} follows from 
\eqref{est:03} due to Sobolev imbedding theorem as $\frac{2}{p}+\frac{3}{q}<1$.
In order to prove \eqref{est:04} we use a fact that
$$
(\vr_1,\ldots \vr_n)=\Psi^{-1}(\vr,\vt_1,\ldots,\vt_{n-1}),
$$
where $\Psi$ is the diffeomorphism defined in \eqref{def:psi},
and therefore 
\begin{equation} \label{5.9} \begin{split}
&\sup_{t \in (0, T)}\|\vr_k(\cdot, t)-\vr_k^0(\cdot)\|_{L_q(\Omega)}
\leq \int^T_0\|\pd_t (\Psi^{-1})({\vec \theta}(\cdot, t), 
\vr_0(\cdot)+\sigma(\cdot, t))\|_{L_q(\Omega)}\,dt\\
&\quad \leq \int^T_0
\|(\Psi^{-1})'({\vec \theta}(\cdot, t), \vr_0(\cdot) + \sigma(\cdot, t))
\|_{L_\infty(\Omega)}
\|(\pd_t{\vec \theta}(\cdot, t), \pd_t\sigma(\cdot, t))\|_{L_q(\Omega)}\,dt. 
\end{split}\end{equation}
By \eqref{est:01} and \eqref{est:04}, we have
\begin{equation} \label{5.8}
\sup_{t \in (0, T)}\|(\Psi^{-1})'({\vec\theta}(\cdot, t), 
\vr_0(\cdot) + \sigma(\cdot, t))\|_{L_\infty(\Omega)} \leq C.
\end{equation}
Thus, by \eqref{5.9} we have
\begin{equation} \label{5.10}\begin{split}
\sup_{t \in (0, T)}\|\vr_k(\cdot, t)-\vr_k^0(\cdot)\|_{L_q(\Omega)} 
&\leq C\int^T_0\|(\pd_t{\vec\theta}(\cdot, t), 
\pd_t\sigma(\cdot, t))\|_{L_q(\Omega)}\,dt \\
& \leq CT^{1/{p'}}\|\pd_t({\vec \theta}, \sigma)\|_{L_p((0, T), L_q(\Omega))}
\leq C(M)E(T).
\end{split}\end{equation}
Similarly, \begin{align*}
&\|\nabla[\vr_k(\cdot, t)-\vr_k^0(\cdot)]\|_{L_q(\Omega)} \\
&\quad
\leq \|(\Psi^{-1})'({\vec \theta}(\cdot, t), \vr^0(\cdot)+\sigma(\cdot, t)
\|_{L_\infty(\Omega)}
\|(\nabla{\vec \theta}(\cdot, t), \nabla\vr^0(\cdot)
+\nabla\sigma(\cdot, t))\|_{L_q(\Omega)}
+ \|\nabla(\vec{\vr^0})\|_{L_q(\Omega)} \leq C(L,M),
\end{align*}
which implies \eqref{est:05}. In order to show \eqref{est:06} note that 
$
W^{\frac{3}{q}+\epsilon}_q(\Omega)\subset L_\infty(\Omega)$  $\forall \epsilon>0,$  
therefore for $\epsilon < 1-\frac{3}{q}$ we have 
\begin{equation}\label{5.12}\begin{split}
&\sup_{t \in (0, T)}\|\vr_k(\cdot, t)
-\vr_k^0(\cdot)\|_{L_\infty(\Omega)}\\
&\quad \leq (\sup_{0 \in (0, T)}
\|\vr_k(\cdot, t)
-\vr_k^0(\cdot)\|_{L_q(\Omega)})^{\theta}\\
&\qquad\times
(\sup_{0 \in (0, T)}
\|\vr_k(\cdot, t)
-\vr_k^0(\cdot)\|_{H^1_q(\Omega)})^{1-\theta} \leq C(M,L)E(T)
\end{split}\end{equation}
with $\theta = 1-(3/q+\epsilon) \in (0, 1)$. This way we obtain \eqref{est:06} and complete the proof.

\qed
\noindent
The next lemma gives bounds on the terms coming from the change of coordinates. 
\begin{lem} \label{l:est_lag}
Let $A_{2\Delta}(\kv)\nabla^2(\cdot),A_{1\Delta}(\kv)\nabla(\cdot),A_{2\dv}(\kv)\nabla^2(\cdot),A_{1\dv}(\kv)\nabla(\cdot)$ be defined in \eqref{a2delta},\eqref{a1delta},\eqref{a2div} and \eqref{a1div}, respectively. Then 
\begin{align} 
&\|A_{2\Delta}\nabla^2\bv,A_{2\dv}\nabla^2\bv\|_{\lpqt}+
\|A_{1\Delta}\nabla \bv,A_{1\dv}\nabla \bv\|_{L_\infty(0,T;L_q(\Omega))}\leq C(M)E(T) \label{est:06}\\
&\|A_{2\Delta}\nabla^2\vt_k,A_{2\dv}\nabla^2\vt_k\|_{\lpqt}+
\|A_{1\Delta}\nabla \vt_k,A_{1\dv}\nabla \vt_k\|_{L_\infty(0,T;L_q(\Omega))}\leq C(M)E(T)   \label{est:07}
\end{align}
for all $k=1,\ldots,n-1$.
\end{lem}
\emph{Proof.}
By \eqref{est:01} and \eqref{a2delta} we have 
$$
\|A_{2\Delta}\nabla^2\bv\|_{\lpqt}\leq \|\bV^0(\kv)\|_{L^\infty(\Omega\times(0,T))}(1+\|\bV^0(\kv)\|_{L^\infty(\Omega\times(0,T))})\|\nabla^2 \bv\|_{\lpqt}\leq C(M)E(T).
$$
Next, notice that  
$$
\left\| \int_0^t \nabla^2 \bv \right\|_{L_q(\Omega)} \leq \int_0^t \|\nabla^2 \bv\|_{L_q(\Omega)} \leq t^{1/p'}\|\nabla^2\bv\|_{\lpqt},
$$
therefore, by \eqref{est:01} and \eqref{est:04}, 
\eqh{
\left\| \nabla_{\kv}V^0_{lm}(\kv) \left[\int_0^t \de_l \nabla \bv ds\right] \frac{\de \bv}{\de y_m} \right\|_{L_q(\Omega)} 
&\leq \|\nabla_{\kv}V^0_{lm}(\kv)\|_{L_\infty(\Omega\times(0,T))} \left\|\int_0^t \nabla^2 \bv\right\|_{\lpqt} \|\nabla \bv\|_{L_\infty(\Omega\times(0,T))} \\
&\leq 
C(M)E(T).
}
The other terms in $A_{1\Delta}\nabla \bv$ have a similar structure, therefore we get 
$$
\|A_{1\Delta}\nabla \bv\|_{L_\infty(0,T;L_q(\Omega))}\leq C(M)E(T).
$$
As $A_{2\dv}(\kv)\nabla^2(\cdot)$ and $A_{1\dv}(\kv)\nabla(\cdot)$ have structure similar to $A_{2\Delta}\nabla^2(\cdot)$ and $A_{1\Delta}\nabla (\cdot)$, respectively, we conclude \eqref{est:06}. Finally, $\vt_k$ have the same regularity as $\bv$ so we obtain $\eqref{est:07}$ in the same way. Proof of Lemma \ref{l:est_lag} is complete.

\qed
\noindent
With these results at hand we can proceed with the proof of Proposition \ref{prop:est}.\\ {\bf Estimate of $f_1(U)$}. Since $f_1(U)$ it is 
exactly the same as in the two species case, we obtain (see \cite{PSZ}) 
\begin{equation} \label{est:1}
\|f_1(U)\|_{\lpqt} \leq C(M,L)E(T).    
\end{equation}
{\bf Estimate of $f_2(U)$}. Let us start with $R_2(U)$ defined in \eqref{lag:7}. 
By \eqref{est:01} we have 
\begin{align*}
\|\bV^0(\kv)\lr{\vr_l-\frac{m_l\vr_l\vr}{\Sigma_\vr}}\nabla \vt_{l-1}\|_{\lpqt} \leq C(M,L)E(T).      
\end{align*}
Applying Lemma \ref{l:est_lag} to the remaining terms we obtain   
$$
\|{\bf R}_2(U)\|_{\lpqt} \leq C(M,L)E(T). 
$$
Next, by \eqref{est:02}
$$
\|\sigma\de_t \bv\|_{\lpqt} \leq \|\sigma\|_{L_\infty(\Omega \times (0,T))} 
\|\de_t\bv \|_{\lpqt} \leq C(M)E(T),
$$
and similarly, using \eqref{est:02}-\eqref{est:05} we get 
$$
\left\| \frac{\sigma}{\Sigma_\vr}\nabla\eta,\sigma_l\nabla\vt_{l-1},\frac{m_l}{\Sigma_\vr}(\vr_l\sigma+\vr^0\sigma_l)\nabla\vt_{l-1},\frac{\vr^0}{\Sigma_\vr}\nabla\vr^0 \right\|_{\lpqt} \leq C(M,L)E(T).
$$
In order to estimate the terms with $\frac{1}{\Sigma_\vr}-\frac{1}{\Sigma_\vr^0}$ we write it as 
$$
\frac{1}{\Sigma_\vr}-\frac{1}{\Sigma_\vr^0} = \frac{\Sigma_\vr^0-\Sigma_\vr}{\Sigma_\vr \Sigma_\vr^0}. 
$$
As the denominator is bounded from below by a positive constant,
using \eqref{est:04} we get
$$
\|\vr^0\nabla\eta \left(\frac{1}{\Sigma_\vr}-\frac{1}{\Sigma_\vr^0} \right)\|_{\lpqt} \leq C \sum_{k=1}^3\|\nabla \eta (\vr_k-\vr^0_k)\|_{\lpqt}\leq
$$$$
C \sum_{k=1}^3 \left[ \int_0^T \|\vr_k-\vr^0_k\|_{L_\infty}^p\|\nabla \eta\|_{L_q}^p \right]^{1/p} \leq \sum_{k=1}^3 \|\vr_k-\vr^0_k\|_{L_\infty(H^1_q)}\|\nabla\eta\|_{\lpqt}
\leq C(M,L)E(T), 
$$
and similarly 
$$
\|m_l\vr^0\vr^0_{l}\left( \frac{1}{\Sigma_\vr}-\frac{1}{\Sigma_\vr^0} \right)\nabla\vt_{l-1}\|_{\lpqt} \leq C(M,L)E(T).
$$
Collecting all above estimates we get 
\begin{equation} \label{est:2}
\|{\bf f}_2(U)\|_{L_p(0,T;L_q(\Omega)^3)} \leq C(M,L)E(T).    
\end{equation}
{\bf Estimate of $f_3(U)$}. First we estimate $R^k_3(U)$ given by \eqref{lag:8}-\eqref{lag:10}. For this purpose we show
\begin{lem}
We have
\begin{align} 
\|\cB_{kl}\|_{L_\infty(\Omega\times(0,T))}\leq C(M), \label{est:3_1} \\
\|\nabla \cB_{kl}\|_{\lpqt} \leq C(M)E(T).   \label{est:3_2}
\end{align}
\end{lem}
\emph{Proof.} \eqref{est:3_1} follows directly from \eqref{est:05} and the form of $\cB_{kl}$ \eqref{lag:5b}.
To show \eqref{est:3_2} we need a bound on $\nabla D_{kl}$. For this purpose notice that, by \eqref{est:05}, 
$$
\|\nabla \vr_k\|_{\lpqt}^p \leq \int_0^T\|(\nabla\vr_k - \nabla \vr^0_k)(t,\cdot)\|_{L_q}^p dt + \int_0^T \|\nabla \vr^0_k\|_{L_q(\Omega)}dt \leq [C(M,L)E(T)]^p.
$$
Therefore, under the assumption \eqref{nablaD} and using the fact that the fractional densities are bounded from below by a positive constant we obtain  
\eqref{est:3_2}.

\qed


\noindent
From \eqref{est:07} and \eqref{est:3_1} we get 
\begin{equation} \label{est:3_3}
\|\cB_{kl}(A_{2\Delta}(\kv)\nabla^2\vt_l+A_{1\Delta}(\kv)\nabla \vt_l)\|_{\lpqt} \leq C(M,L)E(T).
\end{equation}
Next, by \eqref{est:04} and \eqref{est:3_2},
$$
\|\nabla \cB_{kl}\nabla \vt_l\|_{\lpqt} \leq \|\nabla\vt_l\|_{L_\infty(\Omega \times (0,T))}\|\nabla \cB_{kl}\|_{\lpqt}\leq C(M,L)E(T).
$$
Therefore 
\begin{align} \label{est:3_4}
&\|\bV^0(\kv)\nabla \cB_{kl}\left([\nabla\vt_l+\nabla\vt_l\bV^0(\kv)]\right)+(\nabla \cB_{kl})\bV^0(\kv)\nabla\vt_l \|_{\lpqt} \nonumber\\
&\leq C \|\|\nabla \cB_{kl}\nabla \vt_l\|_{\lpqt} \leq C(M,L)E(T).
\end{align}
Combining \eqref{est:3_3} and \eqref{est:3_4} we get 
\begin{equation} \label{est:3_5} 
\|R^{kl}_3(U)\|_{\lpqt} \leq C(M,L)E(T).
\end{equation}
Finally, by \eqref{est:01} and \eqref{est:04},
\begin{align*}
\left\|\left(\vr_{k+1}-\frac{m_{k+1}\vr_{k+1}\vr}{\Sigma_\vr}\right)\sum_{j,m=1}^3V^0_{jm}(\kv)\frac{\de v_j}{\de y_m}\right\|_{\lpqt} & \leq C \sum_{j,m=1}^3\|V^0_{jm}(\kv)\|_{L_\infty(\Omega \times (0,T))}\left\|\frac{\de v_j}{\de y_m}\right\|_{\lpqt}\\ & \leq C(M)E(T),
\end{align*}
which together with \eqref{est:3_5} yields
\begin{equation}
\|R^k_3(U)\|_{\lpqt} \leq C(M,L)E(T).
\end{equation}
The remaining terms in \eqref{lin1:5c} contains only components 
of type $\vr_k \nabla \bv, \vr_k \de_t \theta_l, \nabla \vr_k \nabla \theta_l$ and $\vr_k \nabla^2 \theta_l$, therefore we can estimate them in a similar way to ${\bf f}_2(U)$ using \eqref{est:02}-\eqref{est:06} obtaining   
\begin{equation} \label{est:3}
\|f^k_3(U)\|_{\lpqtn} \leq C(M,L)E(T), \quad k=1,\ldots,n-1. 
\end{equation}

\noindent
{\bf Estimate of $f^k_4(U)$}. This task is more delicate since we have to find a bound on 
$\|f^k_4(U)\|_{H^{1/2}_p(\R;L_q(\Omega))}$. 
However, the structure of boundary condition \eqref{bc:normal1} 
is exactly the same as in the two species case, therefore 
we can repeat the estimate from \cite{PSZ}. For the sake 
of completeness we repeat the idea here.
First we have to extend  $f^k_4(U)$ to whole real line.
For this purpose we apply the extension operator \eqref{def:ext}.
Let us denote 
$$
\bJ[\bv](t)= \bn(x)\bV^0(\kv) \left\{\int^1_0(\nabla\bn)
(y + \tau\int^t_0\bv(y, s)\,ds)\,d\tau
\int^t_0\bv(y, s)\,ds\right\}.
$$ 
Then \eqref{lag:9} can be rewritten as 
$$
R^k_4(U)= - \bJ[\bv] \nabla\vartheta_k. 
$$
Since $\bJ[\bv](0)=0$, we can readily define 
\begin{equation}
\tilde \bJ[u]=e_T(\bJ[u])
\end{equation} 
Next, we also need to extend $\vt_k$. The difference is that it does not vanish at $0$, therefore first we first extend the boundary data to  $\tilde \vt_k^0$ defined on $\R^3$ and define 
\begin{equation} \label{Eh}
E\vt_k=e_T [\vt_k^0-\CT(t)\vt_k^0] + \CT(t)\vt_k^0, 
\end{equation}
where $\CT(t)$ is an exponentially decaying semigroup 
(details can be found in Section 5 of \cite{PSZ}).
The norms of extensions are equivalent with the norms defined on $(0,T)$, therefore we have to estimate 
$\|E\bJ[u] \nabla(E\vartheta_l)\|_{H^{1/2}_p(\R,L_q(\Omega))}$.

For this purpose we apply Lemma \ref{lem:5.1}. 
As $\de \Omega$ is uniformly $C^3$, we can extend the normal vector to $E\vn$ defined on $\R^3$ s.t. $\|E\vn\|_{H^2_\infty(\R^3)}\leq C(\Omega)$.
Then we obtain 
\begin{equation} \label{est:10}
\|EJ[\bv]\|_{L_\infty(0,T;H^1_q(\Omega))} \leq C(M)E(T). 
\end{equation}
and, due to \eqref{ext:2},
\begin{equation} \label{est:11}
\|\de_t EJ[\bv]\|_{L_\infty(0,T;L_q(\Omega))}+\|\de_t EJ[\bv]\|_{L_p(0,T;H^1_q(\Omega))} \leq C\left[\|\bv\|_{L_\infty(0,T;H^1_q(\Omega))}+\|\bv\|_{L_p(0,T;H^2_q(\Omega))}\right]\leq C(M).
\end{equation}
In order to estimate $E\nabla \vt_k$ we apply Lemma \ref{lem:5.2} to obtain 
\begin{equation} \label{est:12}
\|E\nabla \vt_k\|_{H^{1/2}_p(0,T;L_q(\Omega))}+\|E\nabla \vt_k\|_{L_p(H^1_q(\Omega))}\leq C(M,L). 
\end{equation}
Applying Lemma \ref{lem:5.1} with $f=E\bJ[u]$ and $g=\nabla(E\vt_k)$ 
and using \eqref{est:10} - \eqref{est:12} we obtain 
\begin{equation} \label{est:13}
\|R^k_4(U)\|_{L_p(\R,H^1_q(\Omega)) \cap H_p^{1/2}(\R,L_q(\Omega))} \leq E(T)C(M,L).
\end{equation}
Now, combining \eqref{est:1},\eqref{est:2},\eqref{est:3},\eqref{est:13} and \eqref{lin1:5d} we obtain \eqref{est:nonlin}, which completes the proof of Proposition \ref{prop:est}.
\subsection{Fixed point argument}
Theorem \ref{thm:lin2} allows us to define an operator 
$(\sigma,\bv,\vt)={\mathcal S}(\osigma,\ov,\ovt)$ as a solution to system \eqref{lin1:sys} with the right 
hand side $f_1(\bar U), {\bf f}_2(\bar U),f^k_3(\bar U),f^k_4(\bar U)$ where $\bar U=(\osigma,\ov,\ovt)$. 
From the Proposition \ref{prop:est} combined with Theorem \ref{thm:lin2} we easily verify that for any $M>0$
$$
{\mathcal S}:\CH_{T,M} \to \CH_{T,M} 
$$ 
is well defined provided $T>0$ is sufficiently small. It remains to show that ${\mathcal S}$ is a contraction on $\CH_{T,M}$. 
For this purpose we show
\begin{prop} \label{prop:est_dif}
Let $\bar U_1=(\osigma_1,\ov_1,\ovt_1),\bar U_2=(\osigma_2,\ov_2,\ovt_2) \in {\mathcal H}_{T,M}$ for given $T,M>0$, where the initial conditions satify the assumptions of Theorem \ref{thm:main2}.  
Let $f_1(U),f_2(U),f_3^k(U)$ and $f^k_4(U)$ be given by \eqref{lin1:5a}-\eqref{lin1:5d}, 
where $R_1(U),R_2(U),R_3^k(U)$ and $R^k_4(U)$ are defined in \eqref{lag:6},\eqref{lag:7},\eqref{lag:8}-\eqref{lag:10} and \eqref{lag:9}, respectively. 
Then 
\begin{align} \label{est:dif1}
&\|f_1(U_1)-f_1(U_2)\|_{L_p(0,T;W^1_q(\Omega))} + \|{\bbf_2}(U_1)-{\bbf_2}(U_2)\|_{L_p(0,T;L_q(\Omega)^3)}
+\|f_3^k(U_1)-f_3^k(U_2)\|_{\lpqtn} +\nonumber\\ 
&\|f^k_4(U_1)-f^k_4(U_2)\|_{L_p(0,T,H^1_q(\Omega)^{n-1})} 
+ \|f^k_4(U_1)-f^k_4(U_2)\|_{H^{1/2}_p(\R,L_q(\Omega)^{n-1})}  \leq E(L,M,T) [U_1-U_2]_T.    
\end{align}
\end{prop}
\emph{Proof}.
The precise form of the terms on the left hand side of \eqref{est:dif1} is rather 
complicated, however what is essential is that it contains only the terms which are products of either 
$\ov_1-\ov_2$, $\osigma_1-\osigma_2$ or $\ovt_1-\ovt_2$ multiplied by some quantities which are small
for small times. Therefore, following the lines of the proof of Proposition \ref{prop:est} we obtain \eqref{est:dif1}.

\qed

Now we can subtract systems for $U_1$ and $U_2$ to obtain a linear problem for $U_1-U_2$ with the structure 
of the left hand side that same as in \eqref{lag:sys}, zero initial and boundary conditions and left hand side which is estimated in 
\eqref{est:dif1}. Therefore, combining Proposition \ref{prop:est_dif} and Theorem \ref{thm:lin2} we obtain 
\begin{equation}
[{\mathcal S}(U_1)-{\mathcal S}(U_2)]_T \leq E(T) [U_1 - U_2]_T,
\end{equation} 
which implies that for any $M>0$, $\mathcal S$ is a contraction on $\CH_{T,M}$ for sufficiently small $T$.
Therefore, application of the Banach fixed point theorem to $\mathcal S$ completes the proof of Theorem \ref{thm:main2}.

\qed  

\section*{Appendix: Proof of Proposition \ref{thm:main0}}
\subsection*{Derivation of the normal form}
The proof of Proposition \ref{thm:main0} is split into a couple of steps. 
First we derive the normal form of system \eqref{1.1}. By the change of unknowns \eqref{def:psi} we have
\begin{equation} \label{norm:1} 
[\nabla \vr,\nabla h_1 \ldots \nabla h_{n-1}]^T = A [\nabla \vr_1,\ldots \nabla \vr_n]^T
\end{equation} 
with 
\begin{equation} \label{norm:2}
A = \left( \begin{array}{cc}
1 & 1_{1\times (n-1)} \\[5pt]
\left(-\frac{1}{m_1\vr_1}\right)_{(n-1) \times 1} & {\rm diag}\left(\frac{1}{m_2\vr_2},\ldots,\frac{1}{m_n\vr_n}\right) \\
\end{array}\right).
\end{equation}
The matrix $A$ is diagonal except the first row and first column, which also have quite simple structure. 
It is therefore easy to observe that its inverse reads
\begin{equation} \label{norm:4} 
A^{-1} = \left( \begin{array}{cc}
\frac{m_1\vr_1}{\Sigma_\vr} & \left[\left(-\frac{m_1\vr_1 m_k \vr_k}{\Sigma_\vr}\right)_{k=2 \ldots n}\right]_{1 \times (n-1)}\\[5pt]
\left[\left(\frac{m_k\vr_k}{\Sigma_\vr}\right)_{k=2, \ldots, n}\right]_{(n-1)\times 1} & {\cal R}
\end{array} \right),
\end{equation} 
where 
\begin{equation} \label{def:sigma}
\Sigma_\vr=\sum_{k=1}^3 m_k \vr_k
\end{equation}
and ${\cal R}$ is matrix of dimension $n-1$ given by
\begin{equation} \label{def:Rkl}
{\cal R}_{kl}=m_{k+1}\vr_{k+1}\delta_{kl}-\frac{m_{k+1}m_{l+1}\vr_{k+1}\vr_{l+1}}{\Sigma_{\vr}}, \quad k,l=1,\ldots, n-1.
\end{equation}
Therefore, from \eqref{norm:1} we obtain 
\begin{equation} \label{norm:3}
[\nabla \vr_1,\ldots \nabla \vr_n]^T=A^{-1}[\nabla \vr,\nabla h_1 \ldots \nabla h_{n-1}]^T
\end{equation}
and, analogously, for the time derivative 
\begin{equation} \label{norm:3a}
[\de_t \vr_1,\ldots \de_t \vr_n]^T=A^{-1}[\de_t \vr,\nabla h^1 \ldots \de_t h_{n-1}]^T.
\end{equation}
From \eqref{norm:3}, \eqref{norm:3a}, and \eqref{norm:4} we infer
\begin{equation} \label{norm:4}
\de_t \vr_{k+1} + \vu \cdot \nabla \vr_{k+1} = \frac{m_{k+1}\vr_{k+1}}{\Sigma_\vr}(\de_t \vr+\vu \cdot \nabla \vr)
+\sum_{l=1}^{n-1}{\cal R}_{kl}(\de_t h_l+\vu \cdot \nabla h_l), \quad k=1,\ldots,n-1. 
\end{equation}
However, from \eqref{1.1} we have
$$
\de_t \vr + \vu \cdot \nabla \vr = -\vr \div \vu
$$
as well as 
$$
\de_t \vr_k + \vu \cdot \nabla \vr_k = -\vr_k \Div \vu - \Div \vF_k.
$$
Inserting these relations to \eqref{norm:4} we obtain  
\begin{equation} \label{norm:5}
\sum_{l=1}^{n-1} {\cal R}_{kl}(\de_t h_l+\vu \cdot \nabla h_l)+\left(\vr_{k+1} - \frac{m_{k+1}\vr_{k+1}\vr}{\Sigma_\vr}\right)\div \vu = -\div {\bF}_{k+1}.
\end{equation} 
We can further rewrite the rhs of the above equations. For this purpose we observe that 
$$
-\frac{\nabla p_1}{\vr}\left( \frac{1}{\vr_1} \sum_{l=2}^3 \vr_l C_{kl} + \frac{\vr_1}{\vr_1}C_{k1}\right)=
-\frac{\nabla p_1}{\vr_1}\sum_{l=1}^3 Y_l C_{kl}=0
$$
due to \eqref{prop_C}. 
Therefore, denoting 
\begin{equation} \label{def:barm}
\bar m = \frac{\vr}{p}
\end{equation}
we obtain from \eqref{eq:diff1} 
\eq{ \label{norm:6}
-\bF_k &= \frac{1}{p}\sum_{l=1}^3 C_{kl} \nabla p_l\\
&=
\frac{\bar m}{\vr}\left[\sum_{l=1}^3 C_{kl}\nabla p_l - \nabla p_1\left(\frac{1}{\vr_1}\sum_{l=2}^3 C_{kl}\vr_l + C_{k1}\right)\right]\\
&=\frac{\bar m}{\vr}\sum_{l=2}^3 C_{kl} \left(\nabla p_l-\frac{\vr_l}{m_1}\frac{\nabla \vr_1}{\vr_1}\right)\\
&=\frac{\bar m}{\vr}\sum_{l=2}^3 \vr_l C_{kl}\left( \frac{\nabla \vr_l}{m_l \vr_l} - \frac{\nabla \vr_1}{m_1 \vr_1}\right)\\
&=\frac{\bar m}{\vr}\sum_{l=2}^3\vr_k\vr_l D_{kl}\nabla h_{l-1}.
}
Now let us transform the pressure term, from \eqref{norm:3} we have
\eq{ \label{norm:7}
\nabla p &= \sum_{k=1}^3 \frac{\nabla \vr_k}{m_k}\\
&=\frac{1}{m_1}\left( \frac{m_1\vr_1}{\Sigma_\vr}\nabla \vr -\sum_{k=2}^3 \frac{m_1\vr_1m_k\vr_k}{\Sigma_\vr}\nabla h_k\right)\\
&\quad +\sum_{l=2}^3\frac{1}{m_l}\left( \frac{m_l\vr_l}{\Sigma_\vr}\nabla \vr + m_l\left(\vr_l-\frac{m_l\vr_l^2}{\Sigma_\vr}\right)\nabla h_{l-1}-\sum_{\substack{k>1\\k\neq l}} \frac{m_l\vr_lm_k\vr_k}{\Sigma_\vr}\nabla h_k\right)\\
& =\frac{\vr}{\Sigma_\vr}\nabla \vr + \sum_{k=1}^{n-1}A_k\nabla h_k,
}
where we denoted
\begin{equation} \label{norm:8}
A_k = \vr_{k+1}-\frac{1}{\Sigma_\vr}\left[m_{k+1}\vr_{k+1}^2+m_{k+1}\vr_{k+1}\sum_{l\neq k+1}\vr_l\right]=\vr_{k+1}-\frac{m_{k+1}\vr_{k+1}\vr}{\Sigma_\vr}.
\end{equation}
From \eqref{norm:5}-\eqref{norm:8}
we obtain the explicit form of the symmetrized system \eqref{sys:normal}.

Now we have to rewrite the boundary conditions \eqref{bc} for the symmetrized system \eqref{sys:normal}. 
First note that with equation for $\vr_1$ being omitted, the system \eqref{sys:normal} needs to be supplemented only with the boundary conditions for $n-1$ last species densities; due to \eqref{norm:6} we get
\begin{equation} \label{bc:normalD}
\bu=0, \quad \frac{\bar m}{\vr}\sum_{l=2}^3 \vr_k \vr_l D_{kl}\nabla h_{l-1} \cdot \bn = 0, \quad k=2,\ldots,n, \quad\mbox{on}\ (0,T)\times\partial\Omega 
\end{equation}
which is exactly \eqref{bc:normal} and it is a natural boundary conditions in view of the second order term in \eqref{sys:normal}$_3$.

\subsection*{Coercivity properties}
Recall that Lemma \ref{lem:1} gives a positive lower bound on fractional densities. We are now ready to prove the more direct coercivity of ${\cal R}$.   
Below, $\xi=(\xi_1,\ldots \xi_n)$ is a vector of complex numbers, 
$\overline{\xi}=(\overline{\xi_1},\ldots \overline{\xi_n})$ is a vector of their complex conjugates, 
and $\langle\cdot,\cdot\rangle$ is a scalar product in $\C$.

\begin{lem} \label{l:R} Let assumptions of Lemma \ref{lem:1} be satisfied.
Then there exists a constant $C_1>0$  independent of (x,t) such that 
\begin{equation} \label{coerc:R}
\langle{\cal{R}}(x,t)\xi,\overline{\xi}\rangle \geq C_1|\xi|^2. 
\end{equation}
\end{lem}
\emph{Proof.}
Notice first that ${\cal R}_{kk}>0$ for every $k=1,\ldots, n-1$. We rewrite ${\cal R}_{kk}$ as 
$$
{\cal R}_{kk}=\frac{1}{\Sigma_{\vr}}m_{k+1}\vr_{k+1}(\Sigma_{\vr}-m_{k+1}\vr_{k+1})=
\frac{1}{\Sigma_{\vr}}m_{k+1}\vr_{k+1}\sum_{l=1,\ l\neq k+1}^{n}m_l\vr_l.
$$  
Then we have due to symmetry of ${\cal R}$
\eq{
\langle{\cal R}\xi,\overline{\xi}\rangle&=\sum_{k=1}^{n-1}{\cal R}_{kk}|\xi_k|^2
+  \sum_{l=1}^{n-1}\sum_{k<l}{\cal R}_{kl}(\xi_k\overline{\xi_l}+\xi_l\overline{\xi_k})\\
&\geq
\sum_{k=1}^{n-1}{\cal R}_{kk}|\xi_k|^2- \sum_{l=1}^{n-1}\sum_{k<l}|{\cal R}_{kl}|(|\xi_k|^2+|\xi_l|^2)\\
&=
\frac{m_1\vr_1}{\Sigma_{\vr}}\sum_{k=1}^{n-1}m_{k+1}\vr_{k+1}|\xi_k|^2\\
&
\geq \frac{m_1\vr_1}{\Sigma_{\vr}}{\rm min}_{k\neq 1}\{m_k\vr_k\}|\xi|^2,
}
which proves \eqref{coerc:R}. 

\qed
Although   \eqref{prop_D} implies only semi-definitness of $D \geq 0$, the change of unknowns introduced in the previous section and resulting reduction by one row and column enables to deduce ellipticity of the resulting matrix which follows from the properties of $D$. The next lemma shows the coercivity of $\cB$.
\begin{lem} \label{l:B} Assume that one of Conditions 1,2 from Proposition \ref{thm:main0} hold. Then there exists a constant $C_2>0$ independent of (x,t) such that 
\begin{equation} \label{coerc:B}
\langle \cB(x,t)\xi,\overline{\xi}\rangle \geq C_2|\xi|^2 \quad \forall \; (x,t) \in \Omega \times [0,T].   
\end{equation} 
\end{lem}
%
%

\noindent \emph{Proof.}
It is convenient to rewrite the entries of $\cB$ as 
\begin{equation}
\cB_{kl}=\frac{\vr}{p} Y_{k+1}Y_{l+1}\frac{C_{k+1,l+1}}{Y_{k+1}}=\frac{\vr}{p} Y_{l+1} C_{k+1,l+1}.
\end{equation}
Under Condition 1 we therefore have 
\begin{equation}
\cB=\frac{\vr}{p} \left( 
\begin{array}{cccc}
Y_2Z_2 & -Y_2Y_3 & \ldots & - Y_2Y_n\\ 
-Y_3Y_2 & Y_3Z_3 & \ldots & - Y_3Y_n\\ 
\ldots & & &\\
-Y_nY_2 & \ldots & & Y_nZ_n
\end{array}
\right).
\end{equation}
In order to compute ${\rm det} \,\cB$ we transform the matrix with elementary operations. 
First we add $n-1$ first rows  to the last one. Denoting the new matrix by $\cB^1$ we have  
$$
\cB^1_{nn}=Y_nZ_n-Y_n \sum_{j=2}^{n-1}Y_j = Y_nY_1
$$
and for  $k<n$ we have
$$
\cB^1_{nk}=-Y_nY_k+Y_kZ_k-Y_k\sum_{j \neq k,j\geq2}Y_j=Y_kY_1, 
$$
therefore  
\begin{equation}
\cB^1=\frac{\vr}{p} \left( 
\begin{array}{cccc}
Y_2Z_2 & -Y_2Y_3 & \ldots & - Y_2Y_n\\ 
-Y_3Y_2 & Y_3Z_3 & \ldots & - Y_3Y_n\\ 
\ldots & & &\\
Y_1Y_2 & Y_1Y_3 & \ldots & Y_1Y_n
\end{array} 
\right).
\end{equation}
Notice that all entries of the last column contain $Y_n$ and all entries of the last row contain $Y_1$,
therefore 
\begin{equation} \label{B:2}
{\rm det} \,\cB = \left(\frac{\vr}{p}\right)^{n-1} Y_1Y_n 
{\rm det}\underbrace{\left(
\begin{array}{cccc}
Y_2Z_2 & -Y_2Y_3 & \ldots & - Y_2\\ 
-Y_2Y_3 & Y_3Z_3 & \ldots & - Y_3\\ 
\ldots & & &\\
Y_2 & Y_3 & \ldots & 1
\end{array}
\right)}_{\cB^2}.
\end{equation}
Now we can easily diagonalize part of the above matrix. For this purpose 
we add to each $k$-th row, $k=1 \ldots n-1$, the last row multiplied by $Y_{k+1}$. 
Then all the entries except the diagonal becomes zero. Namely, we have 
$$
\cB^2_{k,\cdot} + Y_{k+1}\cB^2_{n,\cdot} = Y_{k+1}\sum_{j=1}^3Y_j {\bf e}_k.
$$
Therefore \eqref{B:2} yields 
\begin{equation}
{\rm det}\, \cB =  \left(\frac{\vr}{p}\right)^{n-1}\prod_{k=1}^3 Y_k \left(\sum_{k=1}^3Y_k \right)^{n-1} \geq C >0,
\end{equation}
since $Y_k(x,t)>C$  for every $k=1,\ldots,n$ uniformly w.r.t. $(x,t)$, due to \eqref{rhoidown}.
Next, denoting 
\begin{equation} \label{minors}
{\rm det} \,\cB_k = \left| \begin{array}{ccc}
\cB_{11} & \ldots & \cB_{1k}  \\
\vdots&\ddots&\vdots\\
 \cB_{k1}&\ldots   & \cB_{kk}  
\end{array} \right|
\end{equation}
we have ${\rm det}\,\cB_{k}>0$. Therefore, all the leading principal minors of matrix $\cB$ are positive and hence we have shown 
\begin{equation} \label{detB:pos}
\cB(x,t)>0, \quad {\rm det}\cB(x,t) \geq C>0 \quad \textrm{uniformly in} \; (x,t).    
\end{equation}
Now from \eqref{detB:pos} it's easy to deduce \eqref{coerc:B}. For this purpose note that the eigenvectors $\zeta_i(x,t)$ of $\cB(x,t)$ form an orthonormal basis of $\R^3$ and $\cB(x,t)$ in this basis is in a form 
\begin{equation} \label{B:diag}
\cB(x,t)={\rm diag}(\lambda_1(x,t),\ldots,\lambda_n(x,t)\},\quad \lambda_i(x,t)\geq C>0 \quad \textrm{uniformly in} \;(x,t).
\end{equation}
Therefore, denoting $\xi=\sum_{i=1}^3 \alpha_i\zeta_i$ we 
have 
$$
\langle B(x,t)\xi,\bar\xi\rangle=\sum_{i=1}^3\lambda_i(x,t)\alpha_i^2 \geq {\rm min}_i \{\lambda_i(x,t)\}\sum_{i=1}^3 \alpha_i^2 \geq C |\xi|^2 \quad \textrm{uniformly in} \; (x,t).
$$

Now let us consider a general form of $ D$ satisfying the assumptions \eqref{prop_D}. In this case we use the form of $\cB$ as in \eqref{lag:5b}. In particular, each entry of $k$-th row of $\cB$ contains $Y_{k+1}$, therefore 
\begin{equation}
{\rm det} \,\cB = \lr{\frac{\vr^2}{p}}^{n-1}Y_2 \ldots Y_n \left|
\begin{array}{cccc}
Y_2D_{22} & Y_3D_{23} & \ldots & Y_n D_{2n}\\
Y_2D_{32} & Y_3D_{33} & \ldots & Y_n D_{3n}\\
\ldots & & &\\
Y_2D_{n2} & Y_3D_{n3} & \ldots & Y_n D_{nn}
\end{array}
\right|
\end{equation}
Similarly, since each entry of $k$-th column contains $Y_{k+1}$, we have
\begin{equation}
{\rm det} \,\cB = \left(\frac{\vr^2}{p}\right)^{n-1} (Y_2 \ldots Y_n)^2 \,
{\rm det} \underbrace{ \left(
\begin{array}{cccc}
D_{22} & D_{23} & \ldots & D_{2n}\\
D_{32} & D_{33} & \ldots & D_{3n}\\
\ldots & & &\\
D_{n2} & D_{n3} & \ldots & D_{nn}
\end{array}
\right)}_{:=\bar D}
\end{equation}
Due to \eqref{rhoidown} we have $Y_2 \ldots Y_n \geq C>0$, and so, the whole coefficient in front of matrix $\bar D$ is positive. Notice however that we only have $D \geq 0$ in general, but $\bar D$ is a $(n-1)\times(n-1)$ sub-matrix of $D$
for which we can show positive definiteness. Assume on the contrary that there is a vector $ [v_2,\ldots,v_n]\neq 0$, s.t. 
\begin{equation*}
\bar D [v_2,\ldots,v_n] = 0. 
\end{equation*} 
Then one would also have  that
$$
 D[0,v_2,\ldots,v_n]=0,
$$
which is in contradiction with the fact that ${\rm Ker}D = {\rm lin}\{\vec{Y}\}$ and all $Y_k$ are strictly positive.
Similarly we show that the minors \eqref{minors} are positive, hence we conclude that
\begin{equation} \label{Dpos}
 D(x)>0.
\end{equation} 
Now, as for each $(x,t)$ fixed, $D(x,t)$ is a linear operator, we have
\begin{equation} \label{coerc:2}
\forall (x,t) \in \Omega \; \exists c(x,t)>0 \; s.t. \; \langle\bar D(x,t)\xi,\bar \xi\rangle \, \geq \, c(x,t)|\xi|^2, 
\end{equation}
where 
$$
c(x,t)={\rm min}_{|\xi|=1}\langle\bar D(x,t)\xi,\bar \xi\rangle.
$$ 
Finally, if Condition 2 is satisfied, we can have the function $c(x,t)>0$ defined on a compact 
set $\overline{\Omega}\times [0,T]$, hence 
$$
\exists \kappa>0: \;c(x,t) \geq \kappa \quad \forall \; (x,t) \in \overline{\Omega} \times [0,T],
$$  
which completes the proof.  

\qed

\begin{rmk}  
The method which we applied for the special structure \eqref{Cform} can be to some extent 
repeated for a general matrix using the fact that ${\rm Ker} D={\rm lin}\{\vec{Y}\}$. However, 
in the last step we do not obtain a diagonal sub-matrix but just a matrix with modified entries. 
For this matrix coercivity probably could be shown under some additional assumptions on $D$ also for unbounded domain, we leave this direction for further investigation in the future.
\end{rmk}


\begin{thebibliography}{10}
\bibitem{Ad} 
R. A. Adams, J. F. Fournier, 
\newblock \emph{Sobolev Spaces}, 
\newblock Second edition. Pure and Applied Mathematics (Amsterdam), 140. Elsevier/Academic Press, Amsterdam, (2003).


\bibitem{Amann}
H. Amann.
\newblock {\em Quasilinear parabolic problems via maximal regularity}. 
\newblock Adv. Differential Equations 10, 1081--1110, 2005.

\bibitem{BdV1}
H. Beirao da Veiga, R. Serapioni, A. Valli.
\newblock {\em On the motion of nonhomogeneous fluids in the presence of diffusion.} 
\newblock J. Math. Anal. Appl. 85 (1982), no. 1, 179-191.

\bibitem{BdV2}
H. Beirao da Veiga.
\newblock
{\em Diffusion on viscous fluids. Existence and asymptotic properties of solutions.} 
\newblock Ann. Scuola Norm. Sup. Pisa Cl. Sci. (4) 10 (1983), no. 2, 341-355.

\bibitem{BdV3}
H. Beirao da Veiga. 
\newblock {\em Long time behaviour of the solutions to the Navier-Stokes equations with diffusion.} 
\newblock Nonlinear Anal. 27 (1996), no. 11, 1229–1239

\bibitem{B2010}
D.~Bothe.
\newblock {\em On the {M}axwell-{S}tefan approach to multicomponent diffusion.}
\newblock In {\em Parabolic problems}, Vol.~80 of {\em Progr. Nonlinear
  Differential Equations Appl.}, pages 81--93. Birkh{\"a}user/Springer Basel
  AG, Basel, 2011.
  
\bibitem{BD2015}  
D. Bothe, W. Dreyer.   
\newblock {\em Continuum thermodynamics of chemically reacting fluid mixtures}.
\newblock  Acta Mech., 226:1757--1805, 2015.
  
\bibitem{BP2017}
D. Bothe, J. Pr\"{u}ss.
\newblock {\em Modeling and analysis of reactive multi-component two-phase flows with mass transfer and phase transition--the isothermal incompressible case}. 
\newblock Discrete Contin. Dyn. Syst. Ser. S 10, no. 4, 673--696, 2017.

\bibitem{BH15}
M. Bulicek, J. Havrda.
\newblock {\em On existence of weak solutions to a model describing compressible mixtures with thermal diffusion cross effects.} 
\newblock Z. Angew. Math. Mech. 95, 589--619, 2015.


\bibitem{CJ13}
X.~Chen, A.~J{\"u}ngel.
\newblock {\em Analysis of an incompressible {N}avier-{S}tokes-{M}axwell-{S}tefan system.}
\newblock Comm. Math. Phys., 340 (2), pp. 471--497, 2015.



\bibitem{DHP} 
R.~Denk, M.~Hieber, J.~Pr\"u\ss.
\newblock {\em $\CR$-boundedness, Fourier multipliers and problems of 
elliptic and parabolic type}. 
\newblock Memoirs of AMS. Vol 166. no. 788,  2003.

\bibitem{DDGG}
W. Dreyer, P-\'E. Druet, P. Gajewski, C. Guhlke.
\newblock {\em Existence of weak solutions for improved Nernst--Planck--Poisson models of compressible reacting electrolytes.}
\newblock preprint WIAS, 2016.



\bibitem{ES1} 
Y.~Enomoto, Y.~Shibata.
\newblock {\em On the ${\mathcal R}$-sectoriality and the initial 
boundary value problem for the viscous compressible fluid flow}.
\newblock Funkcial Ekvac., {56}(3), 441--505, 2013.

\bibitem{FPT}
E.~Feireisl, H.~Petzeltov{{\'a}}, K.~Trivisa.
\newblock {\em Multicomponent reactive flows: global-in-time existence for large data.}
\newblock  Commun. Pure Appl. Anal., 7(5):1017--1047, 2008.

\bibitem{VG0}
V. Giovangigli. 
\newblock {\em Convergent iterative methods for multicomponent diffusion. }
\newblock IMPACT Comput. Sci. Engin. 3, 244--276, 1991.

\bibitem{VG}
V.~Giovangigli.
\newblock {\em Multicomponent flow modeling}.
\newblock Modeling and Simulation in Science, Engineering and Technology.
  Birkh{\"a}user Boston Inc., Boston, MA, 1999.
  
  
\bibitem{GM1} V. Giovangigli, M. Massot.
 \newblock {\em Asymptotic stability of equilibrium states for multicomponent
reactive flows.}
 \newblock Math. Mod. Meth. Appl. Sci. 8, 251--297, 1998.
 
 \bibitem{GM2} V. Giovangigli, M. Massot.
 \newblock {\em The local Cauchy problem for multicomponent reactive flows in
Full Vibrational Nonequilibrium.} 
   \newblock Math. Meth. Appl. Sci. 21, 1415--1439, 1998.
  
  \bibitem{GPZ}
V. Giovangigli, M. Pokorn\'y, E. Zatorska.
\newblock {\em On the steady flow of reactive gaseous mixture.} 
\newblock  Analysis (Berlin) 35, no. 4, 319--341, 2015.


\bibitem{HMPW13}
M.~Herberg, M.~Meyries, J.~Pr{\"u}ss, M.~Wilke.
\newblock {\em Reaction-diffusion systems of Maxwell-Stefan type with reversible
  mass-action kinetics.}
\newblock  Nonlinear Anal. 159, 264--284, 2017.

\bibitem{Jungel}
A. J{{\"u}}ngel.
\newblock {\em  Entropy Methods for Diffusive Partial Differential Equations, SpringerBriefs in Mathematics}.
\newblock Springer 2016.

\bibitem{JS13}
A.~J{{\"u}}ngel, I.V. Stelzer.
\newblock {\em Existence analysis of {M}axwell-{S}tefan systems for multicomponent
  mixtures.}
\newblock  SIAM J. Math. Anal., 45(4):2421--2440, 2013.

\bibitem{K84}
S. Kawashima.
\newblock {\em Systems of Hyperbolic-Parabolic Composite Type, with
Application to the Equations of Magnetohydrodynamics.}
\newblock Doctoral Thesis, Kyoto University, 1984.

\bibitem{KS88} S. Kawashima, Y. Shizuta. 
\newblock {\em On the Normal Form of the Symmetric
Hyperbolic-Parabolic Systems Associated with the Conservation Laws.}
\newblock  Tohoku Math. J., 40, pp. 449--464, 1988.

\bibitem{MT13}
M.~Marion, R.~Temam.
\newblock {\em Global existence for fully nonlinear reaction-diffusion systems describing multicomponent reactive flows.}
\newblock  J. Math. Pures Appl. (9), 104 (1), pp. 102--138, 2015.

\bibitem{MaNi} 
A.~Matusumura, T. Nishida. 
\newblock {\em Initial-boundary value problems for the equations of motion of compressible viscous and heat-conductive fluids.} 
\newblock  Comm. Math. Phys. 89 no. 4, 445--464, 1983.

\bibitem{MPZ}
P.~B. Mucha, M.~Pokorn{{\'y}}, and E.~Zatorska.
\newblock {\em Approximate solutions to model of two-component reactive flow.}
\newblock  Discrete Contin. Dyn. Syst. Ser. S, 7, no.5 , 1079--1099, 2014.

\bibitem{MPZ1}
P.~B. Mucha, M.~Pokorn{{\'y}}, E.~Zatorska.
\newblock {\em Chemically reacting mixtures in terms of degenerated parabolic
  setting.}
\newblock  J. Math. Phys., 54(071501), 2013.

\bibitem{MPZ2}
P.~B. Mucha, M.~Pokorn{{\'y}}, E.~Zatorska.
\newblock {\em Heat-conducting, compressible mixtures with multicomponent diffusion: construction of a weak solution. }
\newblock  SIAM J. Math. Anal. 47, no. 5, 3747--3797, 2015.


\bibitem{Murata} 
M.~Murata.
\newblock {\it On a maximal 
$L_p$-$L_q$ approach to the compressible viscous fluid flow
with slip boundary condition}.
\newblock Nonlinear Analysis, {106}, 86--109, 2014.


\bibitem{MS16}
 M.~Murata,  Y.~Shibata. 
 \newblock {\it On the global well-posedness for the compressible 
Navier-Stokes equations with slip boundary condition.}
\newblock J. Differential Equtions {\bf 260} (7), 5761--5795, 2016.


\bibitem{PP1}
T. Piasecki, M. Pokorn{{\'y}}.
\newblock {\em Weak and variational entropy solutions to the system describing steady flow of a compressible reactive mixture. }
\newblock  Nonlinear Anal. 159, 365--392, 2017.

\bibitem{PP2}
T. Piasecki, M. Pokorn{{\'y}}.
\newblock {\em On steady solutions to a model of chemically reacting heat conducting compressible mixture with slip boundary conditions.}
\newblock{Mathematical analysis in fluid mechanics—selected recent results, 223--242, Contemp. Math., 710, 
Amer. Math. Soc., Providence, RI, 2018.}


\bibitem{PSZ}
T. Piasecki, Y. Shibata, E. Zatorska
\newblock{\em On strong dynamics of compressible two-component mixture flow.}
\newblock{SIAM J. Math. Anal. 51 (4) (2019) 2793--2849.}

\bibitem{PSZ2}
T. Piasecki, Y. Shibata, E. Zatorska
\newblock{\em On the maximal $L_p$-$L_q$  regularity of solutions
to a general linear parabolic system.}
\newblock{preprint: arXiv:1903.11281, 2019}


\bibitem{Pruss}
J. Pr\"uss.
\newblock {\em Maximal regularity for evolution equations in $L_p$-spaces. }
\newblock Conf. Sem. Mat. Univ. Bari 285, 1--39, 2003.

\bibitem{SSZ}
H. Saito, Y.~Shibata, X. Zhang.
\newblock {\em Some free boundary problem for two phase inhomogeneous incompressible flow} 
\newblock{preprint: arXiv:1811.02179, 2018}

\bibitem{S17} 
Y.~Shibata.
 \newblock {\it On the local wellposedness of free boundary problem
for the Navier-Stokes equations in an exterior domain.}
\newblock Communication on Pure and Applied Analysis 17 (4) July (2018), 1681-1721, DOI: 10.3934/cpaa.2018081


\bibitem{SS1} \newblock Y.~Shibata, S.~Shimizu.
\newblock {\it On some free boundary problem for the Navier-Stokes equtions.} 
\newblock Diff. Int. Eqns., {\bf 20}, 241--276, 2007.

\bibitem{SS2} 
Y.~Shibata, S. Shimizu.
{\it On the $L_p$-$L_q$ maximal regularity of the Neumann problem for  the Stokes equations in a bounded domain.} 
J. Reine Angew. Mat., {\bf 615}, 157--209, 2008.

\bibitem{St1}
G.~Str\"ohmer. 
\newblock
{\it About a certain class of parabolic-hyperbolic systems of differential equation.} 
Analysis  {\bf 9}, 1--39, 1989. 

\bibitem{Tanabe} 
H.~Tanabe. 
\newblock {\em Functional analytic methods for partial differential 
equations.} 
\newblock Monographs and textbooks in pure and
applied mathematics, Vol 204, Marchel Dekker, Inc. New York, Basel, 
1997.

\bibitem{EZ}
E.~Zatorska.
\newblock {\em On a steady flow of multicomponent, compressible, chemically reacting
  gas.}
\newblock  Nonlinearity, 24:3267--3278, 2011.

\bibitem{EZ2}
E.~Zatorska.
\newblock {\em On the flow of chemically reacting gaseous mixture.}
\newblock  J. Differential Equations, 253(12):3471--3500, 2012.

\bibitem{EZ3}
E.~Zatorska.
\newblock {\em Mixtures: sequential stability of variational entropy solutions.}
\newblock  J. Math. Fluid Mech. 17, no. 3, 437--461, 2015.

\end{thebibliography}
\end{document}